\documentclass[a4paper]{article}
\usepackage{amssymb}
\usepackage{amsmath}

\usepackage{graphicx}
\usepackage{ifpdf}

\newtheorem{thm}{Theorem}[section]

\newtheorem{lem}[thm]{Lemma}

\newfont{\BBB}{msbm10 scaled \magstep1}
\newfont{\BBS}{msbm10}

\def\a{{\sigma}}
\def\t{{\theta}}

\def\CS{{\mathcal S}}
\def\NN{{\mathbb N}}
\def\RR{{\mathbb R}}
\def\ZZ{{\mathbb Z}}

\newcommand{\wt}{\widetilde}

 \ifpdf
 \DeclareGraphicsExtensions{.pdf,.png,.jpg}
 \else
 \DeclareGraphicsExtensions{.eps}
 \fi
 \graphicspath{{images/}}

\begin{document}

\title{New cubature formulae and hyperinterpolation\\
in three variables
\thanks{Work
supported by the National Science
Foundation under Grant DMS-0604056, by the ``ex-$60\%$'' funds of the
Universities of Padova and Verona,
and by the INdAM-GNCS.}}
\author{\bf{Stefano De Marchi}\\
Dept. of Computer Science, University of Verona (Italy)\\ \\
\bf{Marco Vianello}\\
Dept. of Pure and Applied Mathematics, University of Padova (Italy)\\ \\
\bf {Yuan Xu}\\
Dept. of Mathematics, University of Oregon (Eugene, USA)}
\date{}
\maketitle

\begin{abstract}
A new algebraic cubature formula of degree $2n+1$ for the product 
Chebyshev measure in the $d$-cube with $\approx n^d/2^{d-1}$ nodes
is established. The new formula is then applied to polynomial 
hyperinterpolation of degree $n$ in three variables, in which 
coefficients of the product Chebyshev orthonormal basis are computed 
by a fast algorithm based on the 3-dimensional FFT. 
Moreover, integration of the hyperinterpolant provides a new 
Clenshaw-Curtis type cubature formula in the 3-cube. 
\end{abstract}

\maketitle

\section{Introduction.}

A cubature formula with high accuracy is an important tool for numerical 
computation and has various applications. One of the applications is to
construct polynomial hyperinterpolation,  introduced by 
Sloan \cite{S95}, which is an approximation process 
constructed by applying the cubature formula on the Fourier coefficients of the 
orthogonal projection operator. 

A cubature formula of degree $2n-1$ with $N$ nodes with respect to a 
measure $d\mu$ supported on a set $\Omega$ takes the form
\begin{equation} \label{cub}
\int_\Omega{p(x)\,d\mu}=\sum_{\xi\in X_m}{w_{\xi}\,p(\xi)}
\qquad \mbox{for all} \quad p\in \Pi_{2n-1}^d(\Omega)\;,
\end{equation}
where $\{w_{\xi}\}$, called weights,  are (positive) numbers,  $X_n$ is a 
set of points, called nodes, 
\begin{equation} \label{nodes}
   \xi: = (\xi_1,\xi_2,\ldots,\xi_d) \in X_n \subset \Omega 
\end{equation}
with $\mathrm{card}(X_n) = N$, and $\Pi_m^d$ denoted  the subspace 
of $d$-variate polynomials of total degree $\leq m$ restricted to 
$\Omega$. For a cubature 
formula of degree $2n-1$ to exist, it is necessary that 
\begin{equation} \label{dims}
N:=\mbox{card}(X_n)\geq \mbox{dim}(\Pi_n^d(\Omega))={n+d \choose
d} = \frac{n^d}{d!} (1+ o(1)). 
\end{equation}
There are improved lower bounds of the same order in terms of $n$. A challenging
problem is to construct cubature formulae with fewer nodes, that is, with the
number of nodes $N$ close to the lower bound.  

In this paper we consider the case that the measure is given by the product
Chebyshev 
weight function 
\begin{equation} \label{chebw}
d\mu=W_d(x)\,dx, \qquad W_d(x): =  \frac{1}{\pi^d}
      \prod_{i=1}^d \frac{1}{\sqrt{1-x_i^2}}
\end{equation} 
on the cube $\Omega: = [-1,1]^d$.  For $d =1$, the Gaussian quadrature formula
of degree $2n-1$ needs merely $N =n$ points.  Our main result is a new family
of cubature formulae that uses $N \approx n^d/2^{d-1}$ many nodes. For $d =2$
these formulae are known to have minimal number of nodes. For $d \ge 3$ they 
are still far from the lower bound, but they appear to be the best ones that are
known at this moment. We refer to Section 2 for further discussions. We  present
numerical tests on these cubature formulae in three variables and also apply
them to constructing polynomial hyperinterpolation operator in three variables.   

For every function $f\in C(\Omega)$ the $\mu$-orthogonal projection of $f$ 
on $\Pi_n^d(\Omega)$ is 
\begin{equation} \label{proj}
\mathcal{S}_nf(x)=\sum_{|\alpha|\leq n}{a_\alpha\,p_\alpha(x)}, 
\qquad a_\alpha:=\int_{\Omega}{f(x)\,p_\alpha(x)\,d\mu}\;,
\end{equation}
where $x=(x_1,x_2,\ldots,x_d)$ is a $d$-dimensional point,
$\alpha$ is a $d$-index of length $|\alpha|$
\begin{equation} \label{multiindex}
\alpha=(\alpha_1,\ldots ,\alpha_d)\in \mathbb{N}^d,\qquad
|\alpha|:=\alpha_1+\ldots +\alpha_d\;,
\end{equation}
and the set of polynomials $\{p_\alpha\,,\,0\leq |\alpha|\leq n\}$
is any $\mu$-orthonormal basis of $\Pi_n^d(\Omega)$ with $p_\alpha$ 
of total degree $|\alpha|$ (concerning the theory of multivariate orthogonal
polynomials, we refer the reader to the monograph \cite{DX01}). Clearly, 
$\mathcal{S}_n p=p$ for every $p\in \Pi_n^d(\Omega)$. Given a cubature
formula \eqref{cub} of degree $\leq 2n$, we obtain from (\ref{proj}) the 
polynomial approximation of degree $n$ by the {\em discretized  
Fourier coefficients\/} $\{{c}_\alpha\}$
\begin{equation} \label{hyper}
f(x)\approx \mathcal{L}_nf(x):= \sum_{|\alpha|\leq
n}{{c}_\alpha\,p_\alpha(x)} \;,\;\; {c}_\alpha:= \sum_{\xi\in
X_n}{w_{\xi}\,f(\xi)\,p_\alpha(\xi)}\;,
\end{equation}
where ${c}_\alpha=a_\alpha$ and thus $\mathcal{L}_n p=\mathcal{S}_n
p=p$ for every $p\in \Pi_n^d(\Omega)$. This is the hyperinterpolation 
operator. It satisfies the basic estimate: for every $f\in C(\Omega)$,
\begin{equation} \label{L^2}
\|f-\mathcal{L}_n f\|_{L_{d\mu}^2(\Omega)}\leq 2\sqrt{\mu(\Omega)}
\,E_n(f)\to 0\;,\;\;n\to \infty\;,
\end{equation}
where $E_n(f):=\inf{\{\|f-p\|_\infty\,,\,\,p\in \Pi_n^d(\Omega)\}}$, so that it
converges in mean. The convergence rate can be estimated by a 
multivariate version of Jackson theorem (see, for example, \cite{Ple96}), 
which shows that $E_n(f)=\mathcal{O}(n^{-p})$ for $f\in C^p(\Omega)$, 
$p\in \mathbb{R}^+$. It becomes an effective approximation tool in the 
uniform norm when its operator norm (the so-called  Lebesgue constant)
grows slowly (cf. \cite{R03,SW00,HAC06,CDMV06}). The hyperinterpolation
has been used effectively in several cases: originally for the sphere 
\cite{R03,SW00}, and more recently for the square \cite{CDMMV06,CDMV06}, 
the disk \cite{HAC06}, and the cube \cite{CDMV08}. 
We will use our new cubature formulae to construct a 
hyperinterpolation operator of three variables for the Chebyshev weight
function on the cube. We show that the computation can be carried out
efficiently using the 3-dimensional FFT and that the algorithm can be 
completely
vectorized. We will also present numerical results on hyperinterpolation of 
several test functions.

The paper is organized as follows. In Section 2 we construct new cubature 
formulae and report results of numerical tests, where comparisons are 
made with tensor-product Gauss-Chebyshev formulae.  Hyperinterpolation
in three variables is considered in Section 3, where we show how to 
compute it effectively and report the results of numerical tests. Finally 
in
Section 4, we obtain a new (nontensorial) Clenshaw-Curtis type formula in 
the cube by integrating the hyperinterpolant in Section 3 and show that
it has a clear superiority over tensorial Clenshaw-Curtis and Gauss-Legendre 
cubature on nonentire test integrands, a phenomenon known for 
1-dimensional and 2-dimensional Clenshaw-Curtis formulae (see
\cite{T08,SVZ}). 

\section{Algebraic cubature for the $d$-dimensional Chebyshev
measure.}

We consider cubature formula for the product Chebyshev weight
function (\ref{chebw}), which is normalized so that its integral
over $[-1,1]^d$ is $1$. For $d =1$, we write $w(x) = W_1(x)$.

Let $\Pi_n^d$ denote the space of polynomials of total degree $\leq n$ in
$d$ variables.
We write $\Pi_n$ if $d =1$. The Gaussian quadrature formula for $w$ takes the
form
\begin{equation} \label{Gauss1d}
 \int_{-1}^1 f(x) w(x) dx = \frac{1}{n}
       \sum_{k=1}^n f( \cos \tfrac{ (2k-1) \pi}{2n})\;,
         \qquad \forall f \in \Pi_{2n-1}\;.
\end{equation}
For $d=2$, a cubature formula of degree $2n-1$ needs at least (cf.
\cite{Moller})
\begin{equation} \label{N*}
   N^* = \dim(\Pi_{n-1}^2) + \left \lfloor  \frac{n}{2} \right \rfloor
=\frac{n(n+1)}{2}+\left \lfloor  \frac{n}{2} \right \rfloor
\end{equation}
many nodes. Cubature formulae that attain this lower bound can be
constructed for the product Chebyshev weight $W_2(x)$ (see
\cite{MP,X2} and the references therein) by studying common zeros of
associated orthogonal polynomials. In \cite{BP}, these cubature
rules were derived by an elementary method which depends on a
factorization of the Gauss-Lobatto quadrature into two sums, over
{\it even} indices and {\it odd} indices, respectively. This
factorization method was also used for $d > 2$ in \cite{BP} and
yields a cubature formula of degree $2n-1$ for $W_d$ with roughly
$n^d / 2^{d/2}$ many nodes.

A close inspection of the factorization method shows that it
actually allows us to derive cubature formulae of degree $2n-1$ for
$W_d$ with roughly $2(n/2)^d$ many nodes. This number of nodes is
substantially less than $n^d$ of the product cubature formula or 
$n^d/2^{d/2}$ of the formulae in \cite{BP}, although it likely far from 
optimal as seen from \eqref{dims}.  

We start with the Gauss-Lobatto formula for $w$ on $[-1,1]$. It takes the form
\begin{equation} \label{GaussLobatto}
 \int_{-1}^1 f(x) w(x) dx = \frac{1}{n}
 \left( \frac{1}{2} f(-1)+  \sum_{j=1}^{n-1} f \left( \cos \tfrac{ j
\pi}{n}\right)
  + \frac{1}{2} f(1)\right)   :=I_n f\;,
 \end{equation}
which again holds for all $f \in \Pi_{2n-1}$. We proceed to factor this rule into two
terms. The factorization depends on whether $n$ is even or $n$ is odd.
Define
\begin{align}\label{OandE}
n=2m: \qquad
\begin{split}
   &  I^E_n f := \frac{1}{n} \left(\frac12 f(-1)+
           \sum_{j=1}^{m-1} f \left(\cos \tfrac{2 j \pi}{n}\right) +
           \frac{1}{2} f(1) \right)  \\
           % \quad \hbox{and} \quad
    & I^O_n f: = \frac{1}{n} \sum_{j=1}^{m} f\left(\cos \tfrac{(2j-1)\pi
}{n}\right)
\end{split}
\end{align}
and define
\begin{align}
n=2m-1: \qquad \,
\begin{split}
   &  I^E_n f := \frac{1}{n} \left( \sum_{j=1}^{m-1} f
\left(\cos \tfrac{2j\pi}{n}\right)
       + \frac12 f(1) \right)  \\
           % \quad \hbox{and} \quad
    & I^O_n f: = \frac{1}{n} \left( \frac12 f(-1)+ \sum_{j=1}^{m-1}
          f \left(\cos \tfrac{(2j-1)\pi}{n} \right) \right)\;,
\end{split}
\end{align}
where we use the superscripts $E$ and $O$ to signify that
the sum is taken over even indices or odd indices, respectively.  Evidently,
the quadrature \eqref{GaussLobatto} becomes
$$
    \int_{-1}^1 f(x) w(x) dx =  I^E_n f + I^O_n f\;,
          \qquad \forall f \in \Pi_{2n-1}\;,
$$
by definition.

The Chebyshev polynomials, $T_n$, are orthogonal with respect to $w$ on $[-1,1]$,
$$
   T_n(t) := \cos n \t\;, \qquad\quad t = \cos \t \;.
$$
The following elementary lemma plays a key role in constructing
cubature formulae on $[-1,1]^d$.

\begin{lem} \label{lem1}
For $n \ge 0$ and $k \in \ZZ$,
$$
I^E_n T_k  =  \begin{cases}
          0,   &   k \ne 0 \mod  n  \\
          \frac{1}{2},  & k =0  \mod n
 \end{cases}
\quad\hbox{and}\quad
I^O_n T_k  =  \begin{cases}
          0,   &   k  \ne 0 \mod  n  \\
              \frac{1}{2},  & k =0, 2n , 4n , \ldots \\
              - \frac{1}{2},  & k =0, n , 3n , \ldots  .
 \end{cases}
$$
\end{lem}

\noindent
{\bf Proof.\/}
The proof follows from elementary trigonometric identities. For example,
for $n =2m$, an elementary calculation shows that
$$
   I^O_n T_k = \frac{1}{n} \sum_{j=1}^m \cos k \tfrac{(2j-1)\pi}{2m}
        = \frac{\sin k \pi}{4m \sin \frac{k \pi}{2m}}
         = \frac{\sin k \pi}{2n \sin \frac{k \pi}{n}}
$$
from which $I^O_n T_k =0$ for $k \ne 0 \mod n$ follows immediately.
The case when $k$ is a multiple of $n$ follows from the first equal
sign of the above equation without summing it up. Similarly,
$$
    I^E_n T_k = \frac{1}{n} \left(\frac{1}{2}\cos k \pi +
        \sum_{j=1}^{m-1} \cos k \tfrac{j\pi}{m} + \frac{1}{2}\right)
        = \frac{\sin k \pi  \cos \frac{k \pi}{n}}{2n \sin \frac{k \pi}{n}}\;,
$$
from which the stated result follows. The proof for $n =2m-1$ is similar and
is omitted for brevity.\hspace{0.5cm}{\bf q.e.d.\/}

\vskip0.5cm

Let $\a \in \{E, O\}^d$, that is, $$\a = (\a_1,\ldots,\a_d) \;
\mbox{with}\; \a_i = E \; \mbox{or} \;\a_i =O.$$ For a function $f :
\RR^d \mapsto \RR$, we define the sum
$$
I_n^{\a_1} \cdots  I_n^{\a_d} f
$$
as a $d$-fold multiple sum in which $I^{\a_k}$ is applied to the
$k$-th variable of $f$. Let us define
\begin{equation} \label{sigmaTilde}
\tilde{\a}_i=\left\{\begin{array}{ll} E & \a_i=O \\ O & \a_i=E
\end{array} \right. \;
\end{equation}
For each $\a \in \{E,O\}^d$, we then define
$$
I_{n,d}^\a f : =    I_n^{\a_1} \ldots I_n^{\a_d} f  +
        {I}_n^{\wt{\a}_1} \ldots  {I}_n^{\wt{\a}_d} f\;.
$$
Since the sum introduces a symmetry among $\a \in \{E,O\}^d$, there
are $2^{d-1}$ distinct $I_{n,d}^\a f$ sums.

\begin{thm} \label{MainTh}
For $d \ge 1$ and each $\a \in \{E,O\}^d$, the cubature formula
\begin{equation}\label{cuba-d}
  \int_{[-1,1]^d} f(x) W_d(x) dx = 2^{d-1}I_{n,d}^\a f
\end{equation}
is exact for $f \in \Pi_{2n-1}^d$ and its number of nodes, $N$, satisfies
$$
    N = 2 \left(\left \lfloor \frac{n}{2} \right \rfloor \right)^d
( 1 + o(n^{-1}))\;.
$$
\end{thm}

\noindent
{\bf Proof.\/}
For $k = (k_1, \ldots,k_d) \in \NN_0^d$ let $T_k(x) := T_{k_1} (x_1)
 \cdots T_{k_d}(x_d)$, which is a polynomial of total degree
$|k|:=k_1+\cdots + k_d$. It suffices to establish \eqref{cuba-d} for
$f\in \{T_k: |k| \le 2n-1\}$, since this set is an orthogonal basis of $\Pi_n^d$.
In this case we have
$$
 \int_{[-1,1]^d} T_k(x) W_d(x) dx = 2^{d-1}\left[
     I_n^{\a_1} T_{k_1}\cdots I_n^{\a_d} T_{k_d}+
            I_n^{\wt{\a}_1} T_{k_1}\cdots  I_n^{\wt{\a}_d} T_{k_d}
        \right]\;.
$$
By the orthogonality of $T_k$, the left hand side is equal to 1 if
$k = (0,\ldots,0)$ and zero if $k \ne (0,\ldots,0)$. From the
definition of $I^E_n$ and $I^O_n$, it is evident that $I^E_n 1 =
I^O_n 1 = 1/2$. Hence, for $k =(0,\ldots,0)$, the right hand side is
equal to $2^{d-1} (2^{-d} + 2^{-d}) =1$, verifying the equation for
$k =(0,\ldots,0)$.

Assume now $0< |k| \le 2n-1$. If one of $k_i \ne 0 \mod n$, then
$I_{n,d}^{\a} T_k=0$ by Lemma \ref{lem1}.  We are left with the case
that $k_i =0 \mod n$ for all $i$. Since $|k| \le  2n-1$, there can
be at most one $k_i =n$. Furthermore, $|k| > 0$ shows that there is
exactly one $k_i = n$. Thus the right hand side becomes $I_n^{\a_i}
T_n + I_n^{\wt{\a}_i} T_n = I^E_n T_n + I^O_n T_n$, which is zero as
$I^E_n T_n =1/2$ and $I^O_n T_n = -1/2$ according to the Lemma
\ref{lem1}.\hspace{0.5cm}{\bf q.e.d.\/}

\vskip0.5cm For the case of $d =2$, Theorem \ref{MainTh} contains
two distinct cubature formulae for $\a = (E,E), \, (E,O)$,
respectively, whose number of nodes are either equal to $N^*$ in
\eqref{N*} or $N^* + 1$, those are the ones that have appeared in
\cite{MP, X2}, and later in \cite{BP}, as mentioned earlier. For $d =3$,
there are 4 distinct formulae for $\a = (E,E,E),\, (E,E,O),\,(E,O,E),\,(O,E,E)$, 
respectively. For $n=2m$, the number of nodes is
$$
 N = \frac{(n+1)^3 + (n+1)}{4} % \quad\hbox{or}\quad
$$
for $\a = (E,E,E)$ and
$$
 N = \frac{(n+1)^3 - (n+1)}{4}
$$
for $\a = (E,E,O),\,(E,O,E),\,(O,E,E)$, respectively. 

In order to demonstrate the effectiveness of the new cubature formula, 
we present in Figures 1-2 numerical results of (\ref{cuba-d}) with $\sigma=(E,E,E)$ 
on the integrals of six text functions with respect to the product Chebyshev 
measure on the 3-cube. The first three functions are analytic entire (a polynomial, 
an exponential and a gaussian), whereas the other three are less smooth:
one analytic but not entire (a 3-dimensional version of the Runge test function), 
one $C^\infty$ nonanalytic, and one $C^2$. These functions are analogues of
test functions for algebraic cubature in dimension 1 and 2, see \cite{T08,SVZ}.
We compare them with two natural choices for cubature on a tensor product 
domain: the tensor-product Gauss-Chebyshev and Gauss-Chebyshev-Lobatto
formulae. The results, obtained with Matlab (cf. \cite{G04}), demonstrate 
the superiority of the new formula in all 
cases, especially for the less smooth functions, in terms of number of function
evaluations. It should be pointed out that, however, the superiority 
for 
the less smooth functions arises for {\em even\/} $n$ (a sort of parity 
phenomenon). Other numerical tests (not reported for brevity) have shown  
that 
the cubature formula has the same behavior for $\a = 
(E,E,O),\,(E,O,E),\,(O,E,E)$. 

A natural question associated with cubature formulae  is polynomial 
interpolation.
Let $X_{n-1}$ denote the set of the nodes of the cubature formula 
(\ref{cuba-d}). The 
interpolation
problem looks for a polynomial subspace, $\CS$, of the lowest degree such that
$$
    P (x) = f(x), \qquad x \in X_{n-1}, \quad \forall f \in C(\RR^d)
$$
has a unique solution in $\CS$. In the case of $d =2$, this problem
is completely solved in \cite{X2}, where $\CS$ is a subspace of
$\Pi_n^2$ which includes $\Pi_{n-1}^2$, and compact formulae of the
fundamental interpolation polynomials are also given there. For $d >
2$, however, the problem is much harder, since the number of nodes
of our cubature is far from $\dim(\Pi_n^d)$. For example, if $d =3$,
then $\dim(\Pi_{n-1}^d) = n(n+1)(n+2)/6\approx n^3/6$, whereas our
cubature has $\approx n^3/ 4$ many nodes. The problem essentially
comes down to study the polynomial ideal that has $X_{n-1}$ as its
variety (see \cite{X3}).

A simpler approach to polynomial approximation via these new nodes 
is given by hyperinterpolation, as described in the Introduction. In the next
section we shall apply such a method in the 3-dimensional case.  

\begin{figure}
\centering
\includegraphics[scale=0.60,clip]{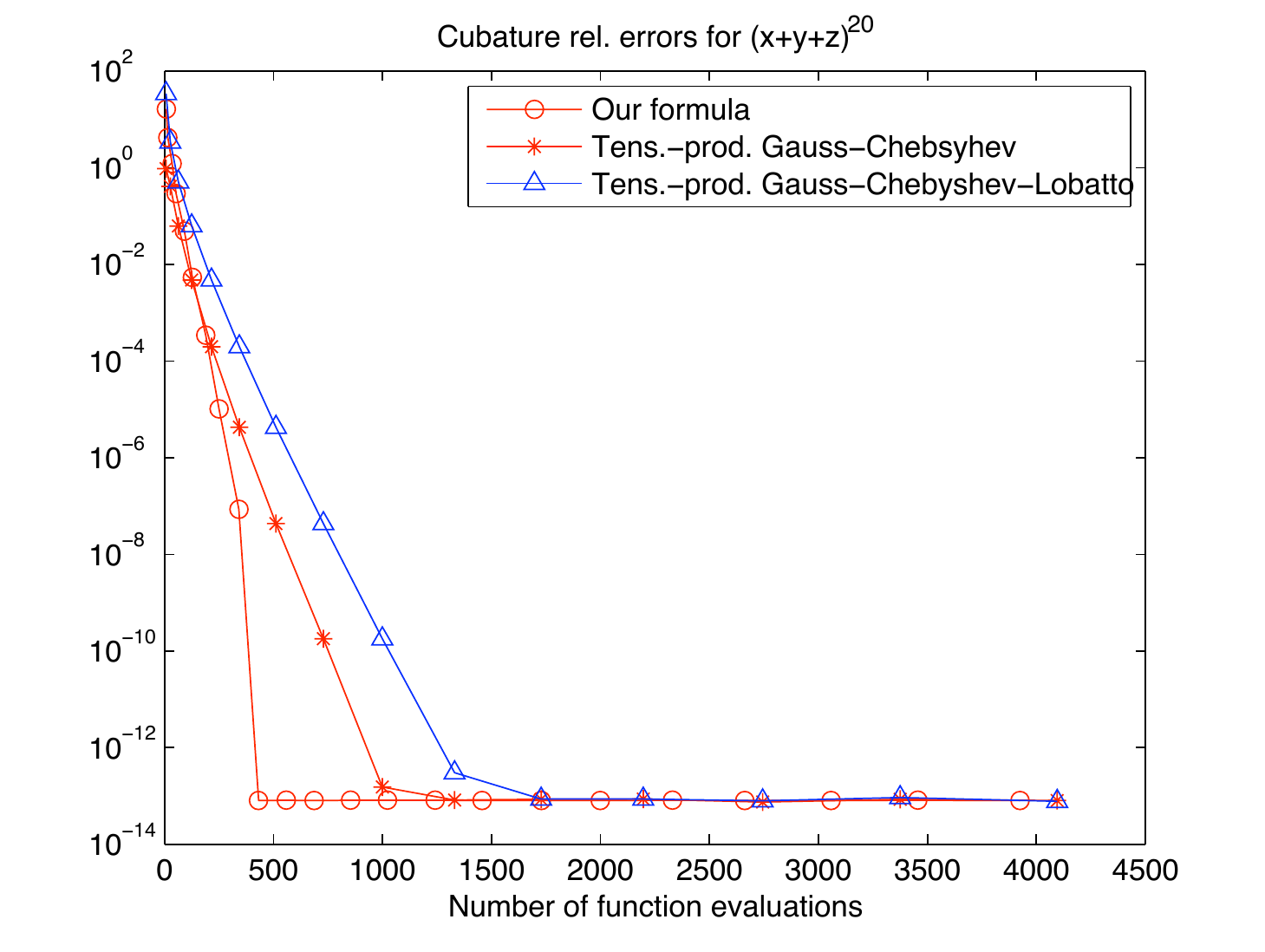}\hfill
\includegraphics[scale=0.60,clip]{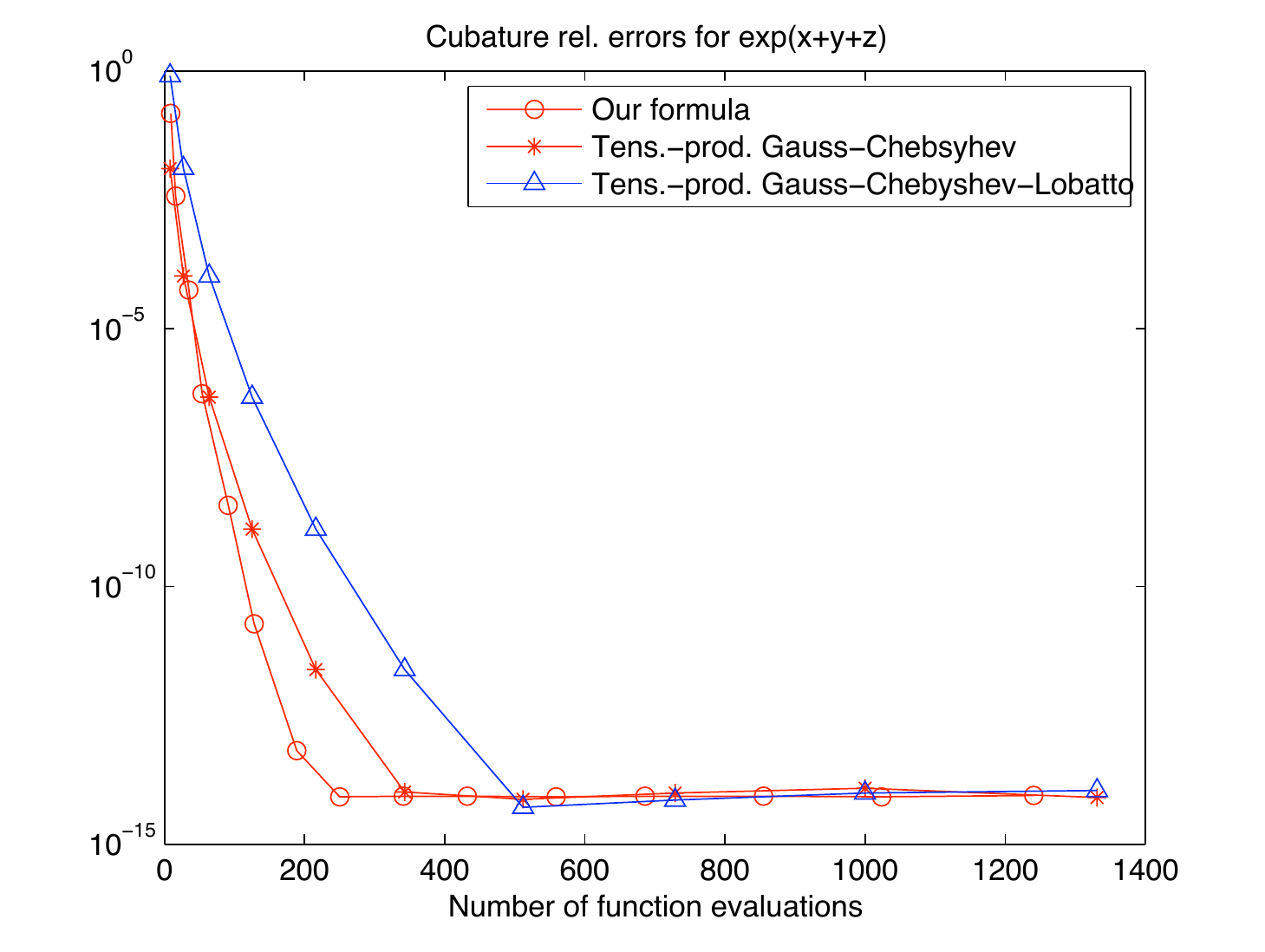}\hfill
\includegraphics[scale=0.60,clip]{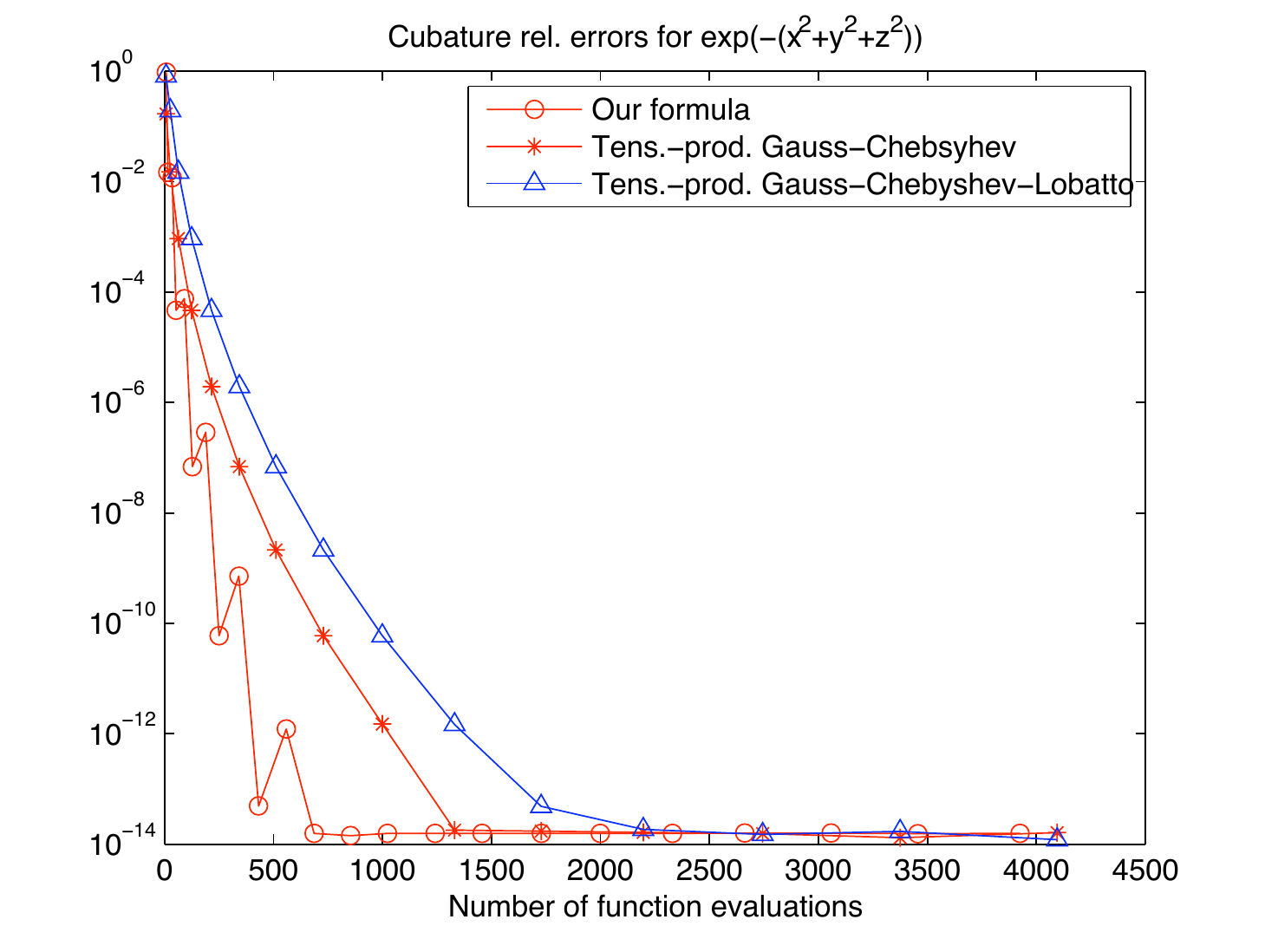}\hfill
\caption{Relative cubature errors versus the number of
function evaluations for three test functions.}
\label{cubchebfig1}
\end{figure}

\begin{figure}
\centering
\includegraphics[scale=0.60,clip]{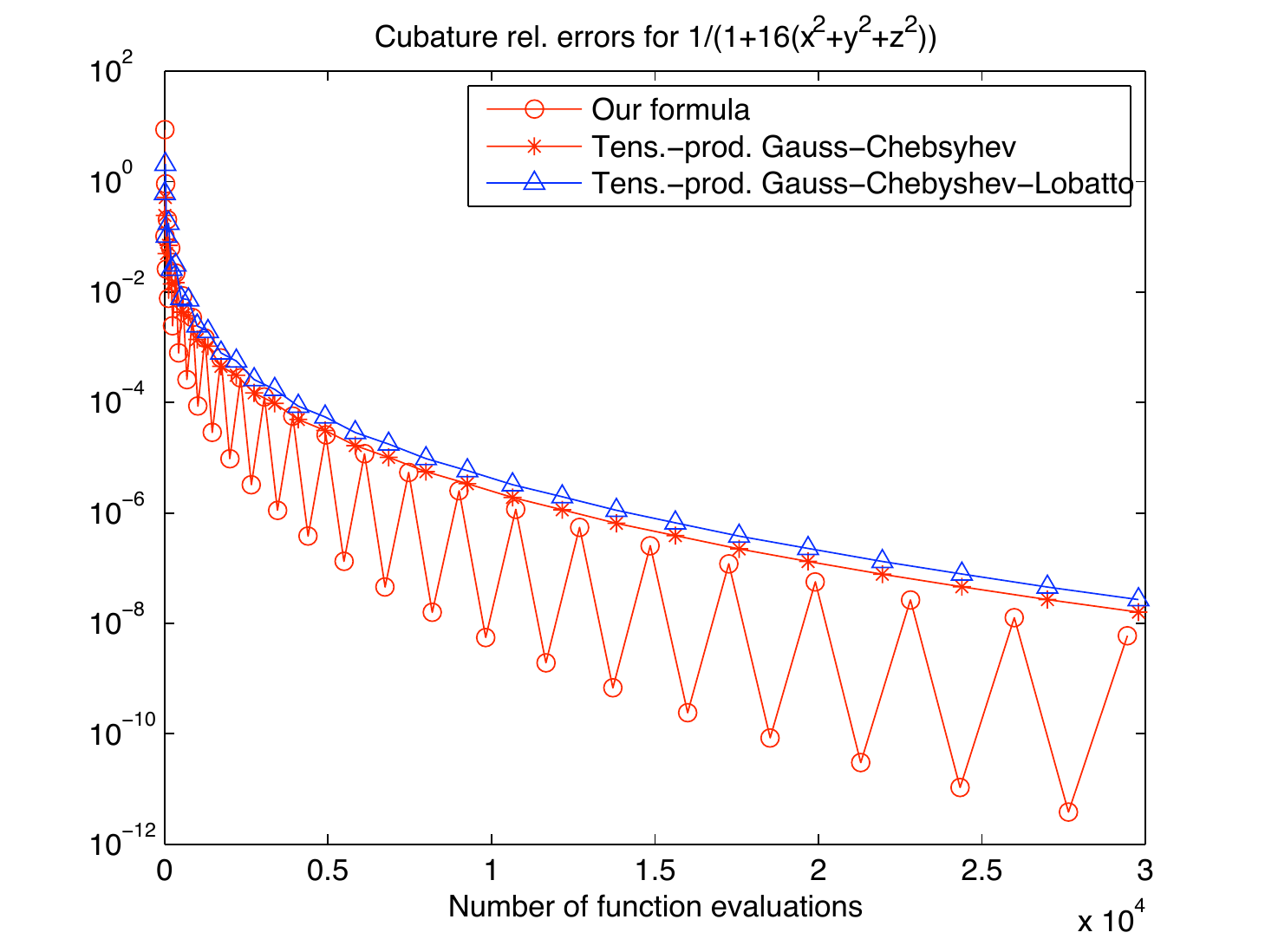}\hfill
\includegraphics[scale=0.60,clip]{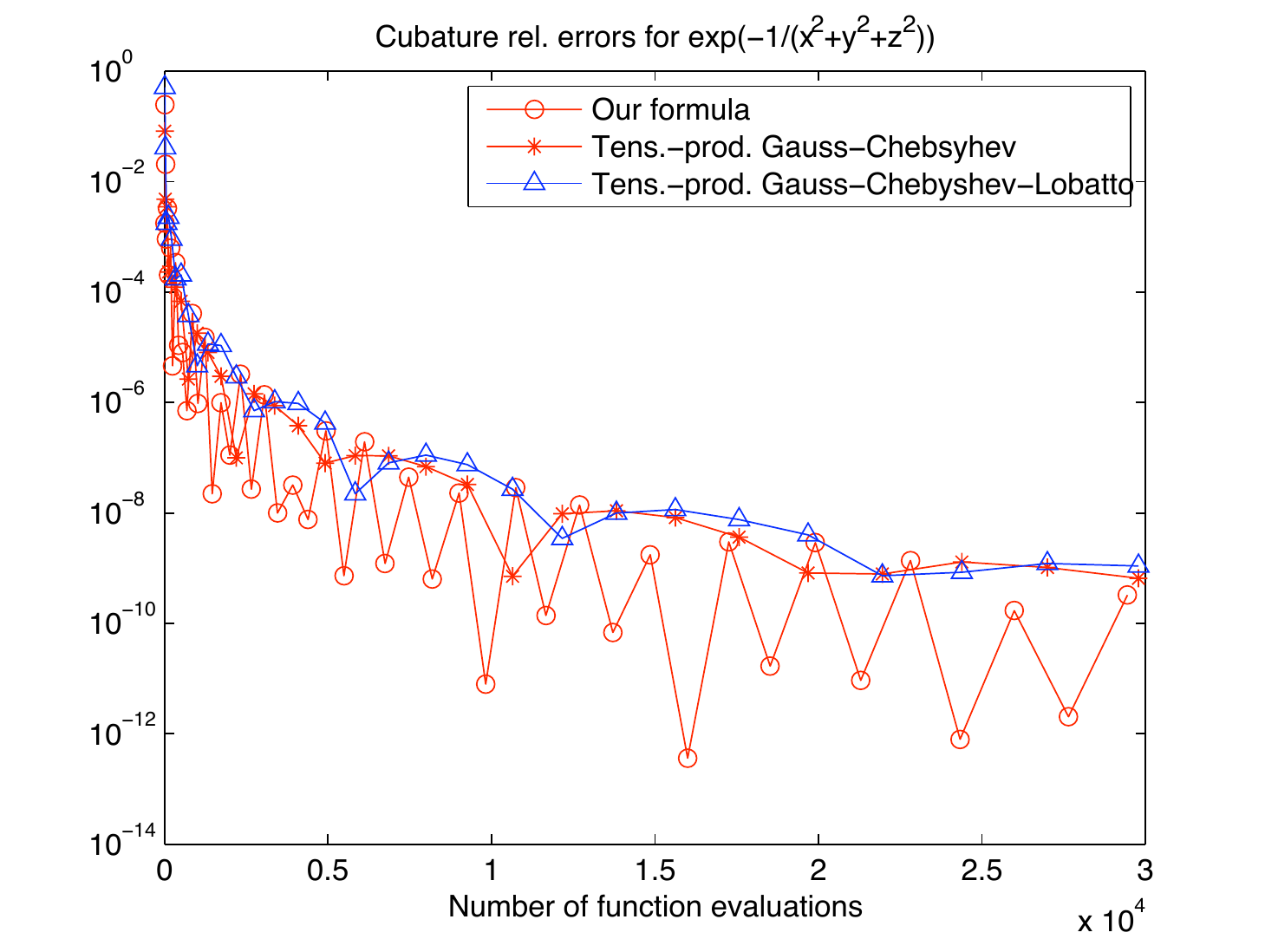}\hfill
\includegraphics[scale=0.60,clip]{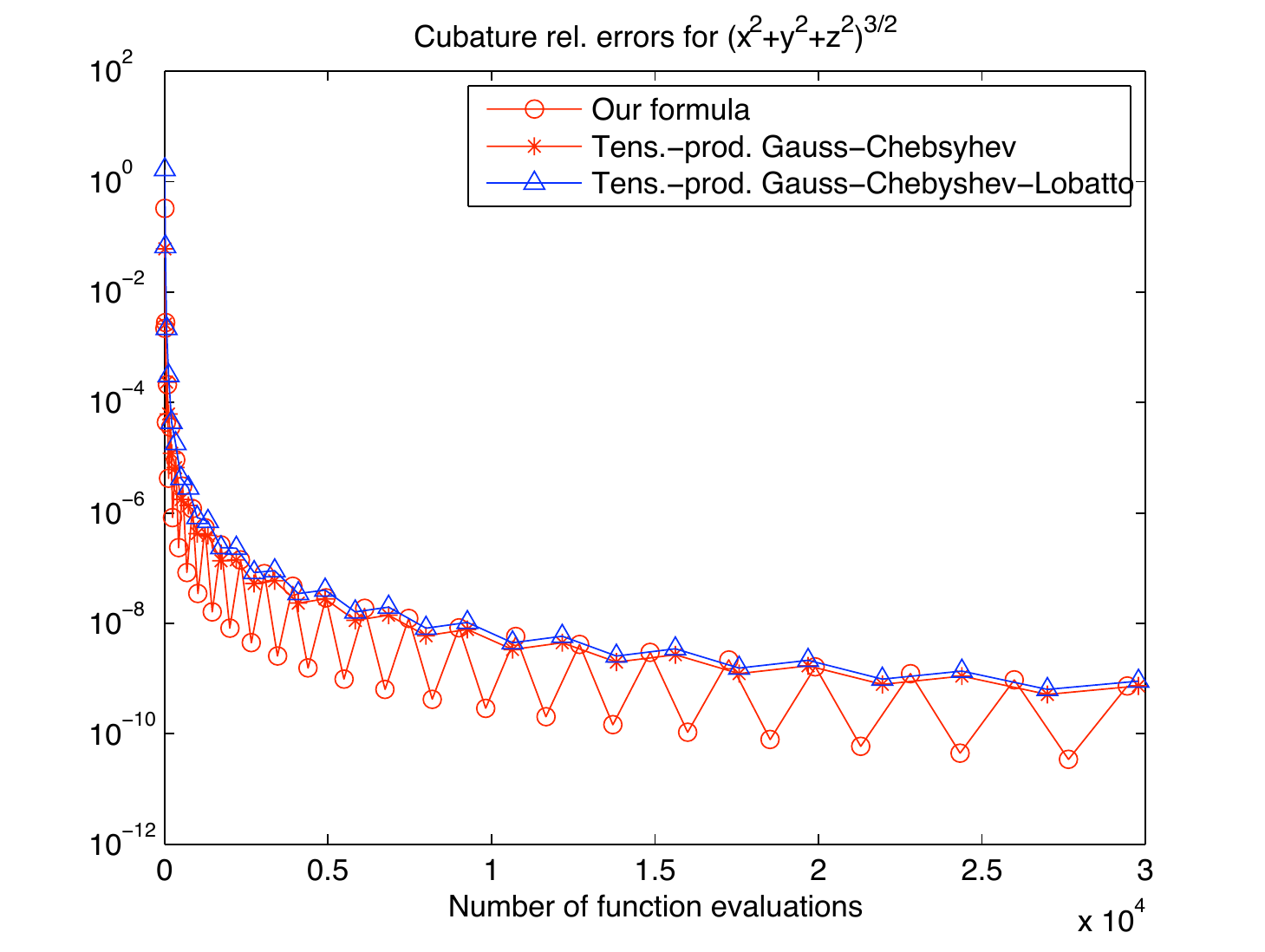}\hfill
\caption{Relative cubature errors versus the number of
function evaluations for three test functions.}
\label{cubchebfig2}
\end{figure}

\section{Implementing hyperinterpolation in the 3-cube.}
We  now use cubature formula (\ref{cuba-d}) to construct hyperinterpolation 
as in (\ref{hyper}) for the 3-cube $\Omega=[-1,1]^3$. In this case, $\{p_\alpha\}$
 is the product Chebyshev orthonormal basis (cf. \cite{DX01}), i.e.
\begin{equation}
p_\alpha(x)=\hat{T}_{\alpha_1}(x_1)\hat{T}_{\alpha_2}(x_2)\hat{T}_{\alpha_3}(x_3)\;,
\end{equation}
where $\hat{T}_{k}(\cdot)=\sqrt{2} \cos(k\arccos(\cdot)), \; k>0$
and $\hat{T}_{0}(\cdot)=1$. Moreover, let
$$C_n=\left\{\cos {k\pi \over n},\;k=0,...,n\right \}\;$$
be the set of $n+1$ Chebyshev-Lobatto points, and $C_n^E$, $C_n^O$
its restriction to even and odd indices, respectively. Then,
\begin{equation} \label{hypernodes}
X_n=\left(C_{n+1}^{\a_1} \times C_{n+1}^{\a_2} \times C_{n+1}^{\a_3} \right)\;
\cup \; \left(C_{n+1}^{\tilde{\a}_1} \times C_{n+1}^{\tilde{\a}_2} \times
C_{n+1}^{\tilde{\a}_3}\; \right)\,,
\end{equation}
with $(\a_1,\a_2,\a_3)\in \{E,O\}^3$, see (\ref{sigmaTilde}). The weights of the 
cubature formula
(\ref{cuba-d}) for $\xi \in X_n$, are
\begin{equation} \label{hyperweights}
w_\xi=\frac{4}{(n+1)^3} \cdot \left\{ \begin{array}{ll}  1 & \mbox{if}\;
\xi \;\mbox{is an interior point} \\
{1 / 2} & \mbox{if}\; \xi\; \mbox{is a face point} \\
{1 / 4} & \mbox{if}\; \xi \; \mbox{is an edge point} \\
{1 / 8} & \mbox{if} \; \xi \; \mbox{is a vertex point} \\
\end{array} \right.
\end{equation}
Note that, since $$\mbox{dim}(\Pi_n^3(\Omega))={(n+1)(n+2)(n+3)/6} < 
N=\mbox{card}(X_n) \approx n^3/4\;,$$ the polynomial
$\mathcal{L}_nf$ in (\ref{hyper}) is not interpolant.

Now, defining
\begin{equation} \label{F} 
F(\xi)=F(\xi_1,\xi_2,\xi_3)=\left\{\begin{array}{ll} w_\xi f(\xi) &
\xi \in X_n \\ \\0 & \xi \in (C_{n+1}\times C_{n+1} \times C_{n+1})
\backslash X_n
\end{array} \right.
\end{equation}
we can write 
\begin{eqnarray*}
c_\alpha&=&\sum_{\xi \in X_n} w_\xi f(\xi) p_\alpha(\xi) \\
& = & \sum_{\xi_1 \in C_{n+1}}\left( \sum_{\xi_2 \in C_{n+1}} \left(
\sum_{\xi_3 \in C_{n+1}} F(\xi_1,\xi_2,\xi_3)\,
\hat{T}_{\alpha_1}(\xi_1)\right) \hat{T}_{\alpha_2}(\xi_2)\right)
\hat{T}_{\alpha_3}(\xi_3) \\
& = & \left(\prod_{s=1}^3{\beta_{\alpha_s}}\right)\,\sum_{i 
=0}^{n+1}\left( \sum_{j 
=0}^{n+1} \left 
(\sum_{k
=0}^{n+1} F_{ijk} \cos{k\alpha_1\pi  \over n+1} \right )
\cos{j\alpha_2\pi  \over n+1} \right) \cos{i\alpha_3\pi  \over
n+1}\;,
\end{eqnarray*}
where 
$$
\alpha=(\alpha_1,\alpha_2,\alpha_3)\in
\{0,1,\dots,n\}^3\;,\;\;
\beta_{\alpha_s}=\left\{\begin{array}{ll}  \sqrt{2} & \alpha_s>0\\
1 & \alpha_s=0
\end{array} \right.\;,\;\;s=1,2,3\;. 
$$
This shows that the 3-dimensional coefficients array $\{c_\alpha\}$ 
is a scaled Discrete 
Cosine Tranform of the 3-dimensional array
\begin{equation} \label{Fijk}
F_{ijk}=F\left(\cos{i\pi  \over n+1},\cos{j\pi  \over n+1},
\cos{k\pi  \over n+1} \right)\;,\;\;0\leq i,j,k\leq n\;,
\end{equation}
where we eventually pick up only the $(n+1)(n+2)(n+3)/6\approx 
n^3/6$ hyperinterpolation coefficients corresponding to 
$|\alpha|=\alpha_1+\alpha_2+\alpha_3\leq n$.  

A fast implementation of hyperinterpolation is now feasible 
(for example in Matlab), via the FFT. 
Indeed, we have written a Matlab code (see \cite{DMV}), completely 
vectorized by several implementation tricks, 
whose kernel can be summarized as follows: 
\vskip0.2cm
\noindent
{\bf Algorithm: Fast total degree hyperinterpolation in the 3-cube\/}
\begin{itemize}
\item[$(i)$] construct the hyperinterpolation point set $X_n$ as union of 
the two 
subgrids in (\ref{hypernodes});

\item[$(ii)$] compute the cubature weights in (\ref{hyperweights});

\item[$(iii)$] compute the 3-dimensional array $\{F_{ijk}\}$ at the 
complete 
grid 
$C_{n+1}\times C_{n+1}\times C_{n+1}$ by (\ref{F}) (notice that $f$ is 
evaluated {\em only\/} at $X_n$); 

\item[$(iv)$] compute the 3-dimensional array of coefficients 
$\{c_\alpha\}$ 
by three nested applications of the 
1-dimensional $\mbox{Real}(\mbox{FFT}(\cdot))$ operator; 

\item[$(v)$] select the coefficients $\{c_\alpha\}$ corresponding to the 
triples $\alpha=(\alpha_1,\alpha_2,\alpha_3)$ such that 
$|\alpha|=\alpha_1+\alpha_2+\alpha_3\leq n$.

\end{itemize}

We recall that there is a simple way to approximate a function in the 
3-cube by  tensor-product of polynomials of degree $n$, that is, by a 
tensor-product discrete Chebyshev series (ultimately a tensor-product 
hyperinterpolant).  Such an approximation uses $(n+1)^3$ function 
evaluations, and $(n+1)^3$ coefficients. In contrast, let us stress 
again the following facts on our total-degree hyperinterpolation of degree
 $n$ in the 3-cube:
\vskip0.2cm
\noindent
{\bf Remark\/}
\begin{itemize}
\item the number of hyperinterpolation nodes, or function 
evaluations, is equal to $\mbox{card}(X_n)\approx n^3/4$;
\item the number of hyperinterpolation coefficients 
is $\mbox{dim}(\Pi_n^3)\approx n^3/6$. 
\end{itemize}

In order to compare the performances of total-degree and tensor-product hyperinterpolation in the 3-cube, we show,  in the following figures, the 
hyperinterpolation errors versus both the number of nodes and the number of coefficients on the six test functions already used in Section 2, and we choose
again $(\a_1,\a_2,\a_3)=(E,E,E)$, see (\ref{hypernodes}). The errors are 
relative to the maximum deviation of the function from its mean and are 
computed on a uniform control grid.  Since the computation of 
the coefficients via the FFT has roughly the same cost for both kinds 
of hyperinterpolation, we have chosen the number of function 
evaluations as a measure of computational cost for the construction, and 
the number of coefficients as a measure of the compression capability of 
the algorithms.  

The situation here is in some sense opposite to that of Figures 1-2. 
Indeed, total-degree appears superior to tensor-product hyperinterpolation 
on the smoothest functions, but not on the less smooth ones. 
As it is natural from the observation above, the behavior of total-degree 
hyperinterpolation in terms of number of coefficients is better than that 
in terms of number of nodes (function evaluations). 

\begin{figure}
\centering
\includegraphics[scale=0.60,clip]{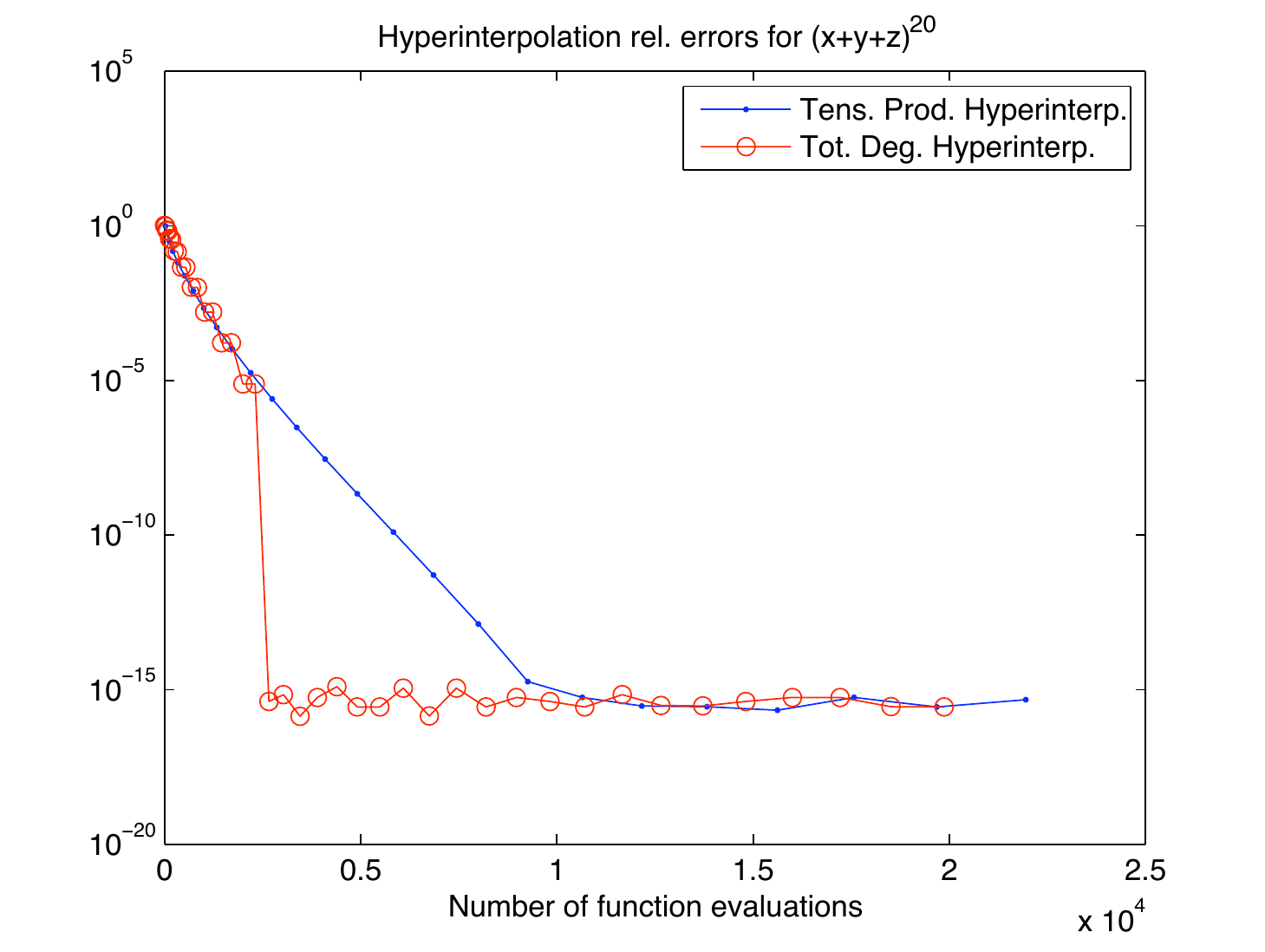}\hfill
\includegraphics[scale=0.60,clip]{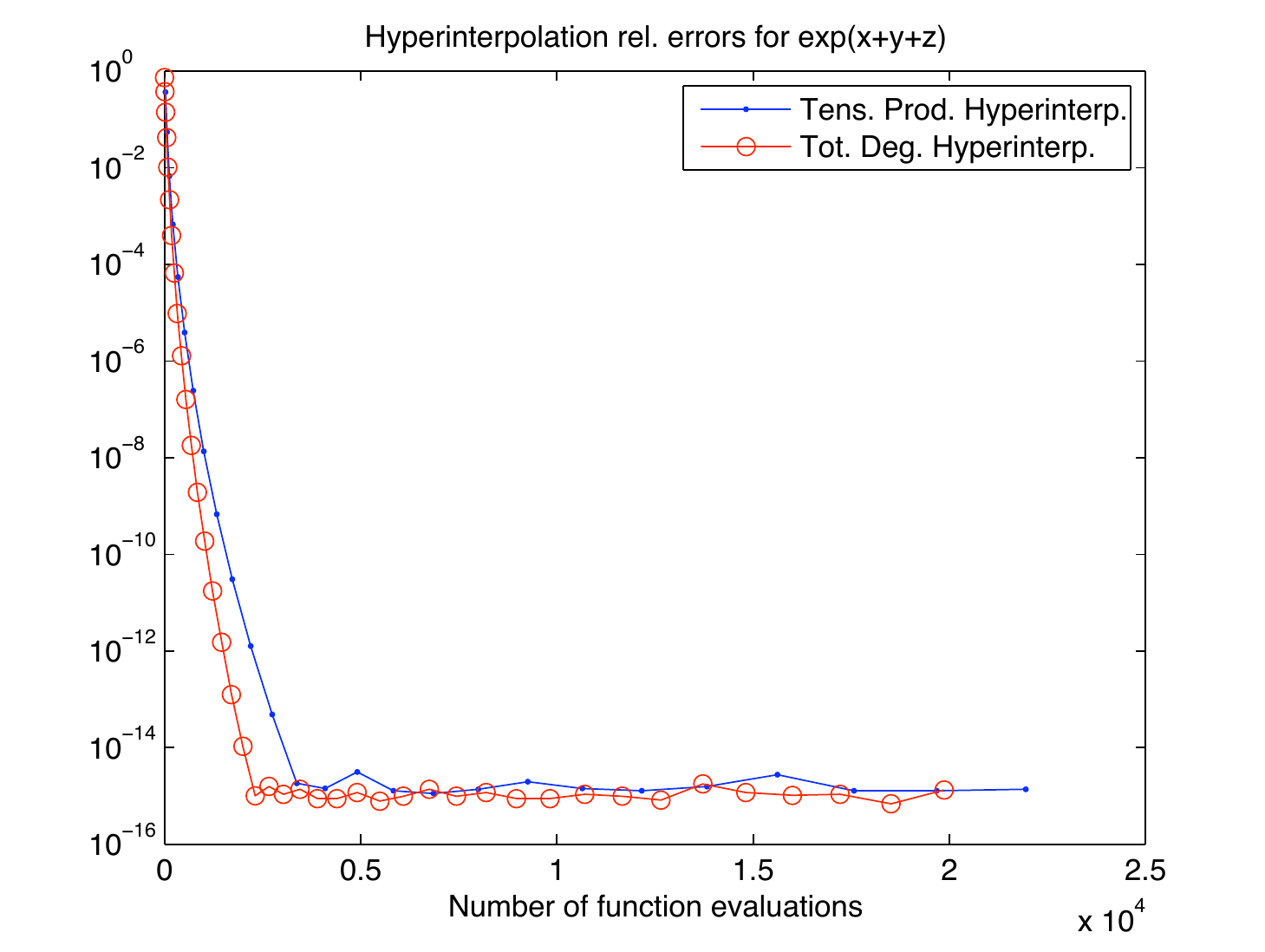}\hfill
\includegraphics[scale=0.60,clip]{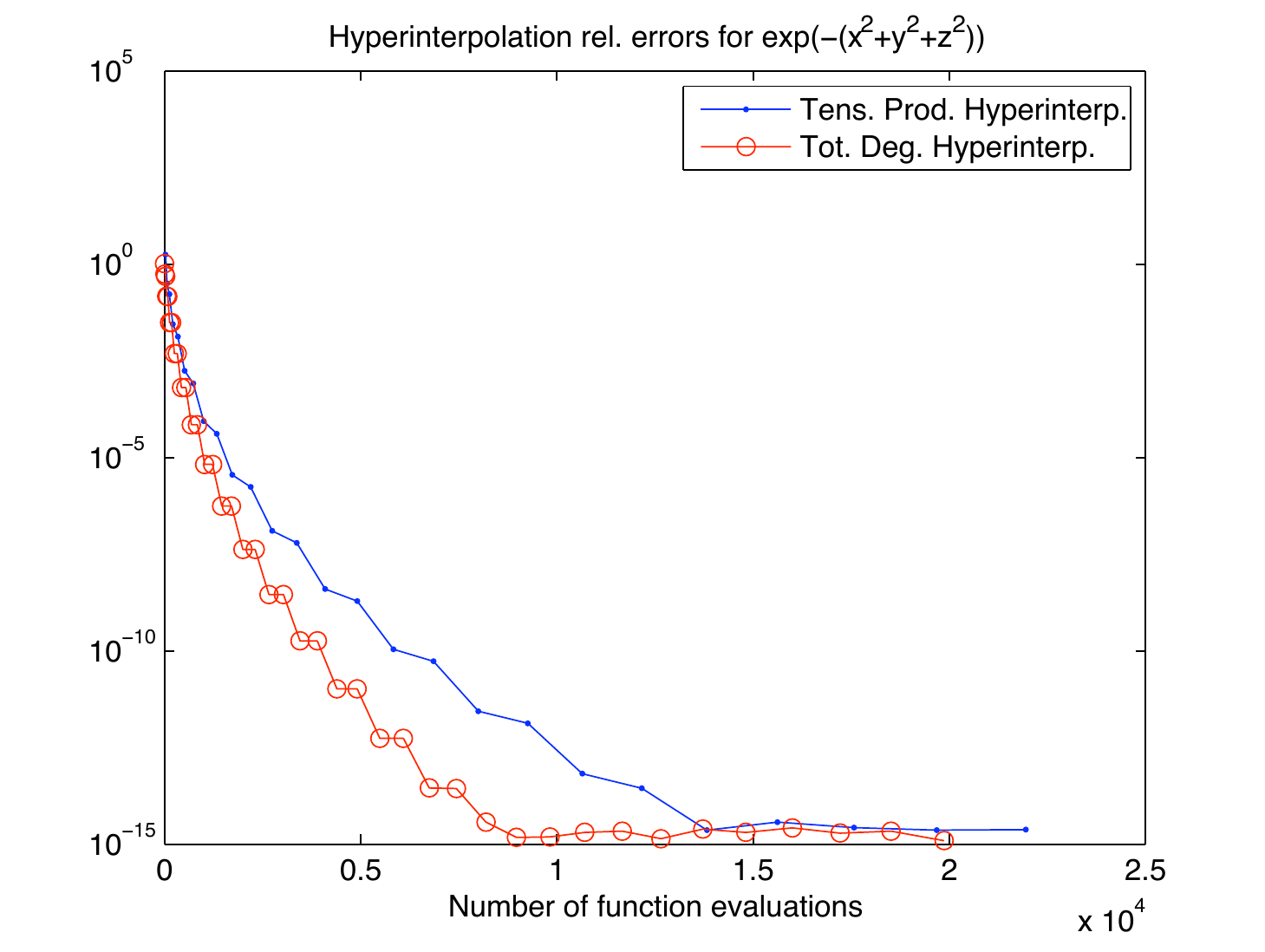}\hfill
\caption{Hyperinterpolation relative errors versus the number of
function evaluations for three entire test functions.}
\label{hyper1fig}
\end{figure}

\begin{figure}
\centering
\includegraphics[scale=0.60,clip]{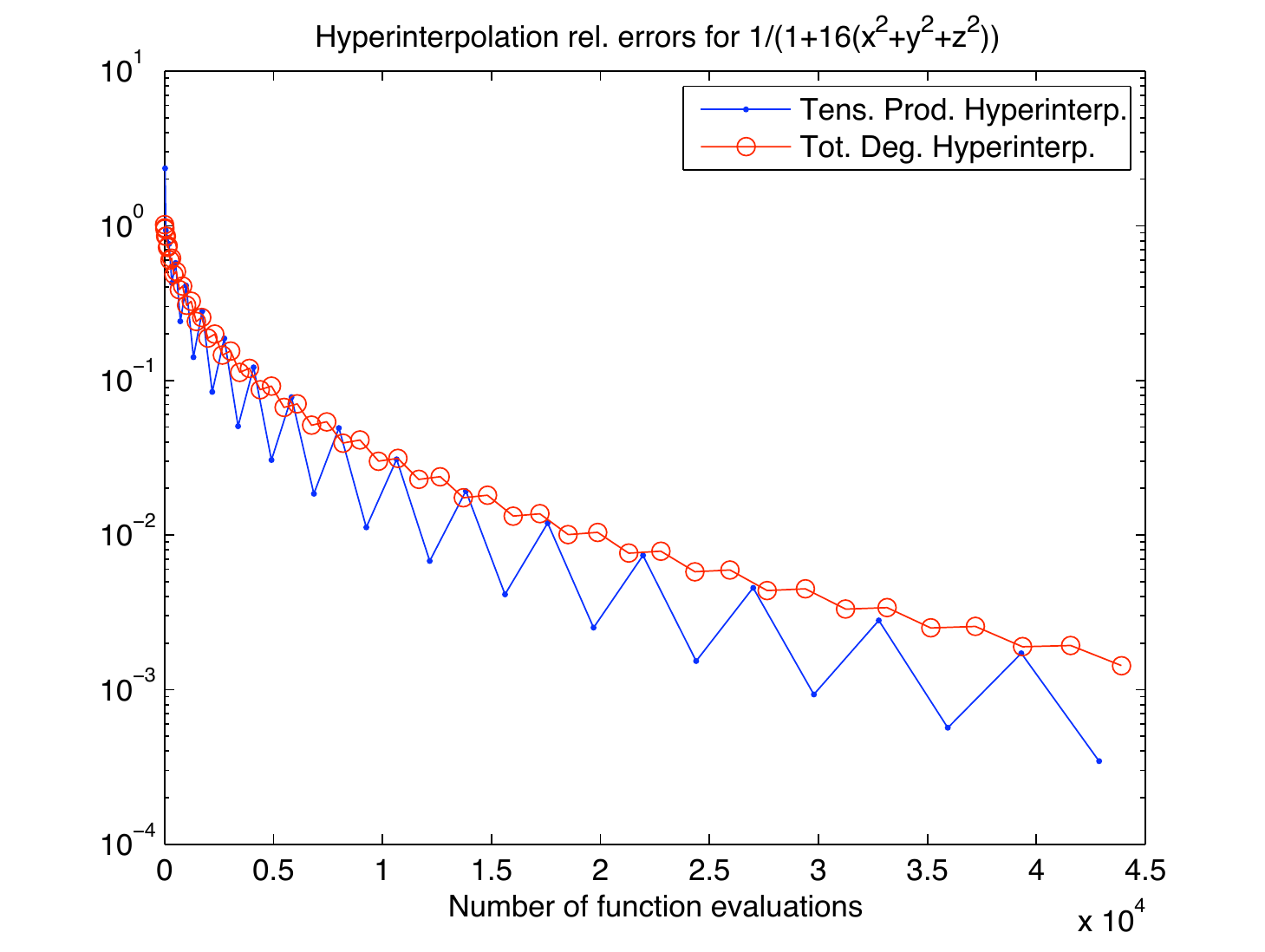}\hfill
\includegraphics[scale=0.60,clip]{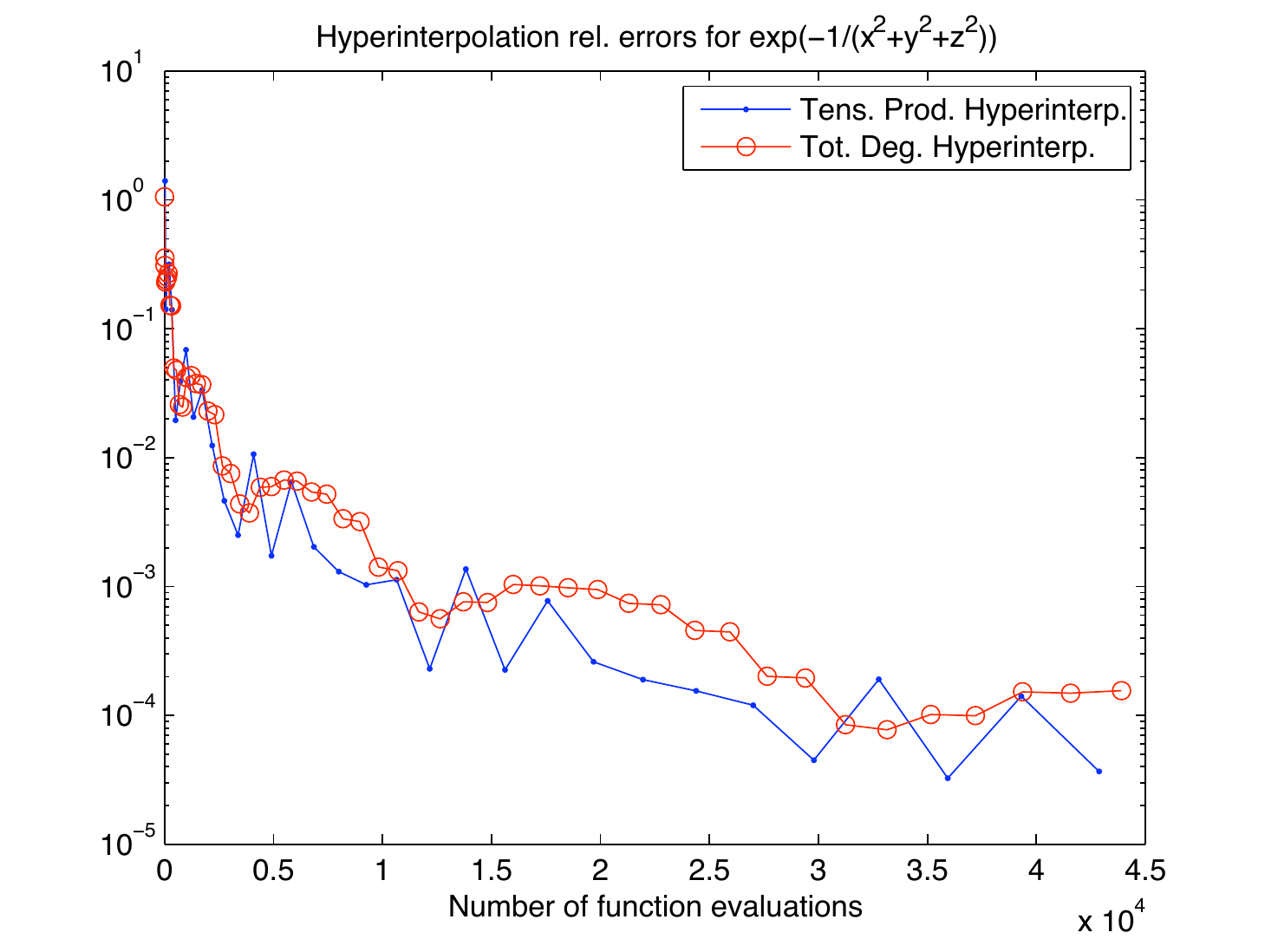}\hfill
\includegraphics[scale=0.60,clip]{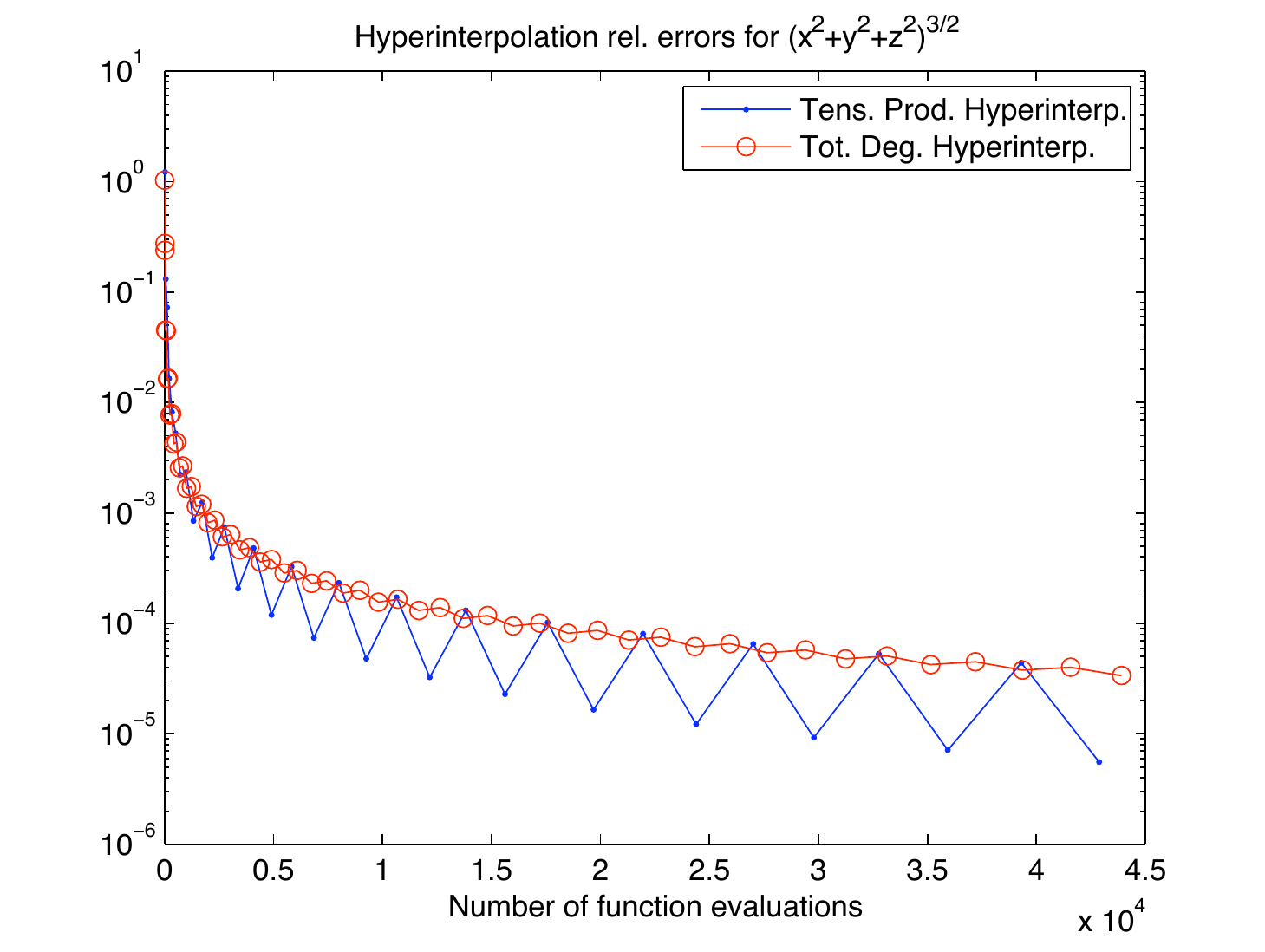}\hfill
\caption{Hyperinterpolation relative errors versus the number of
function evaluations for three nonentire test functions.}
\label{hyper2fig}
\end{figure}

\begin{figure}
\centering
\includegraphics[scale=0.60,clip]{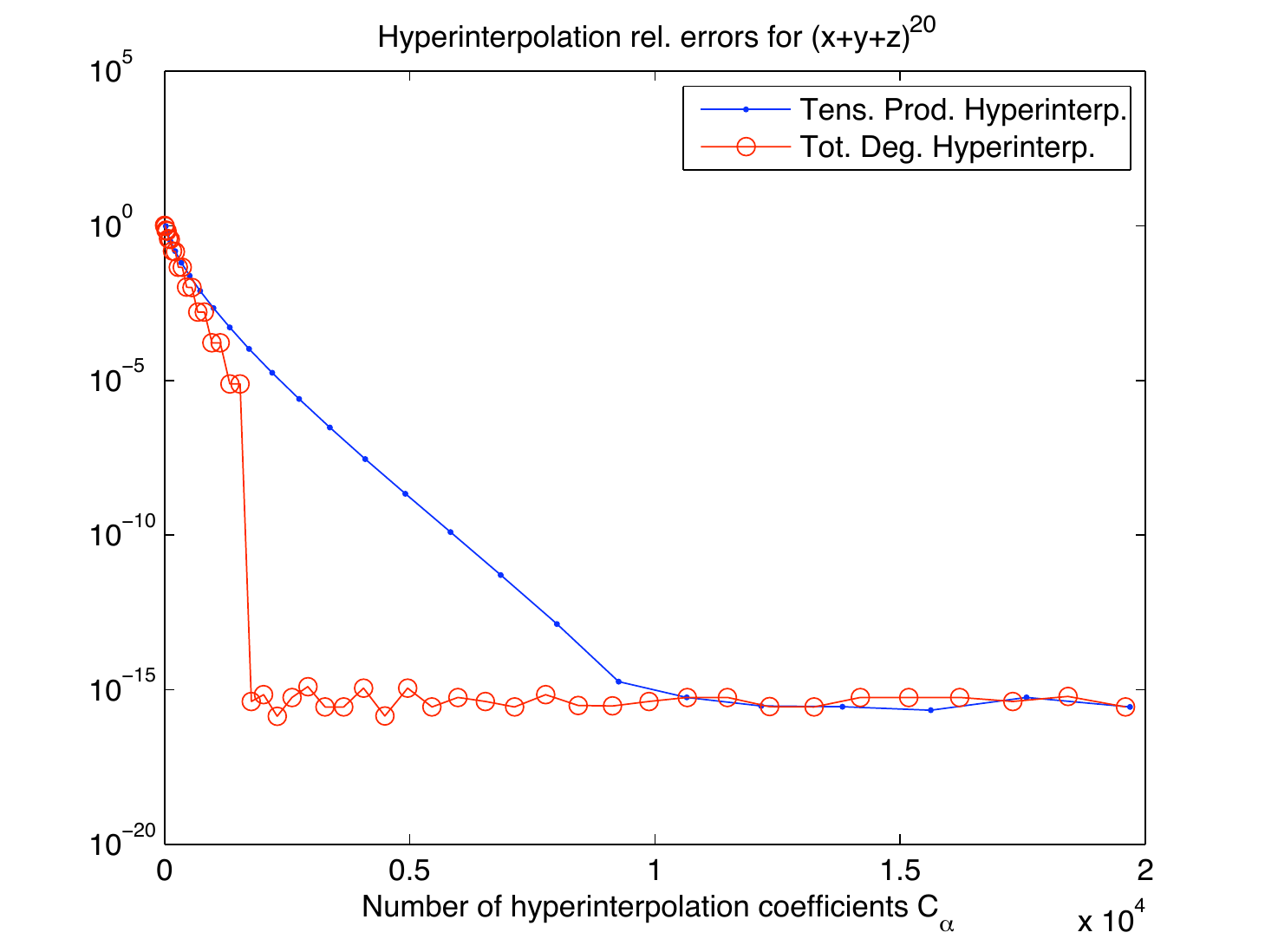}\hfill
\includegraphics[scale=0.60,clip]{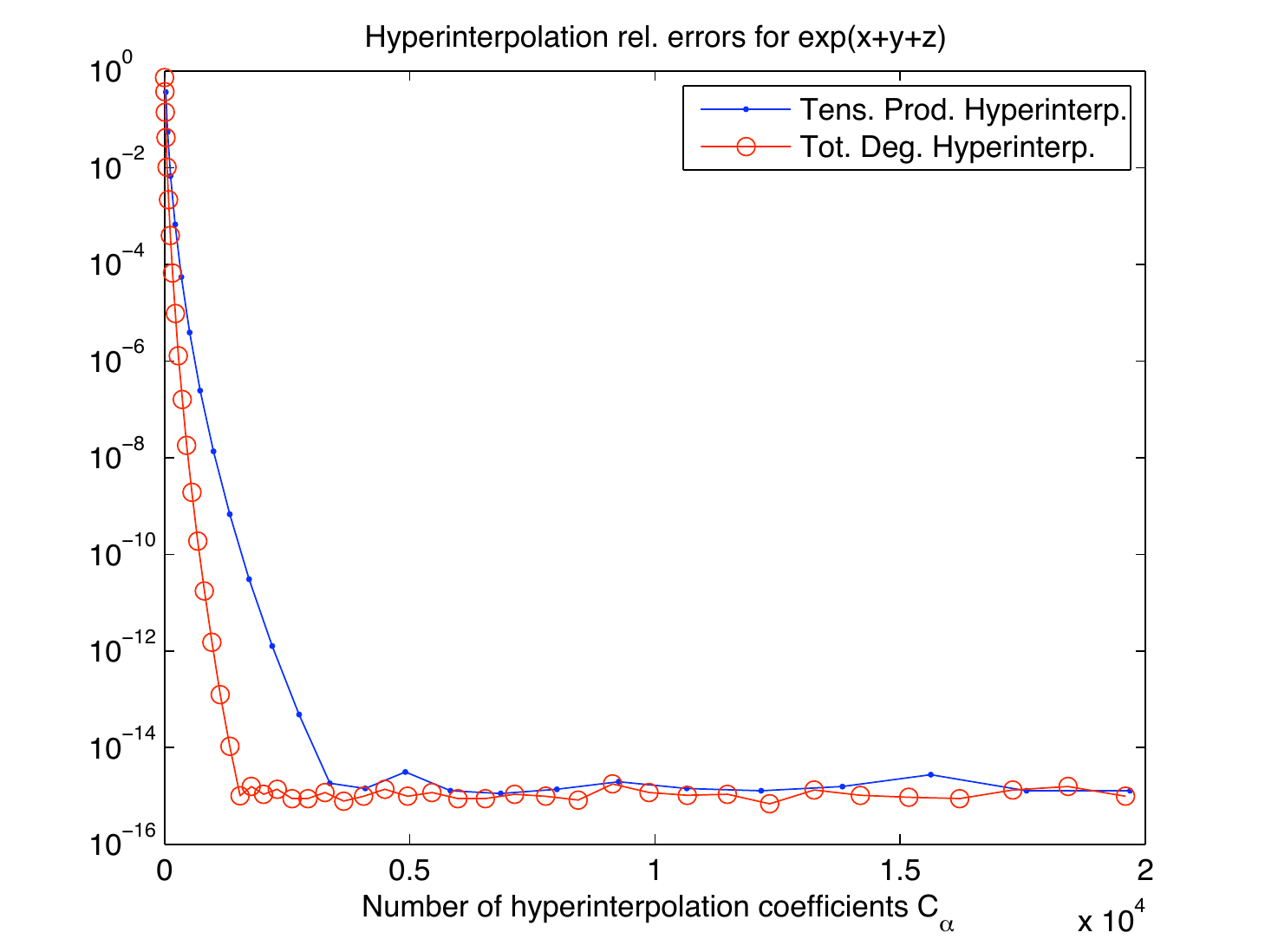}\hfill
\includegraphics[scale=0.60,clip]{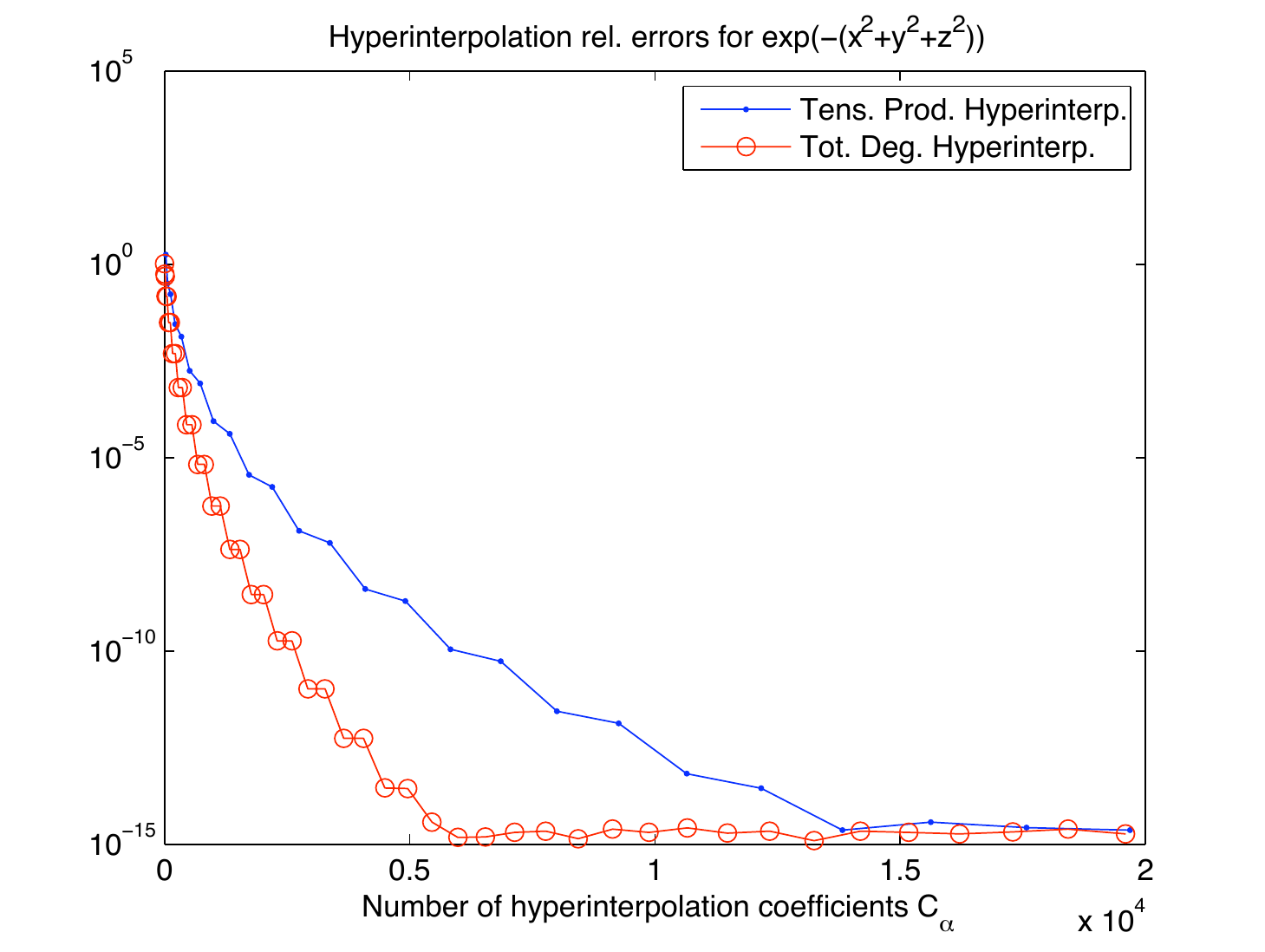}\hfill
\caption{Hyperinterpolation relative errors versus the number of
hyperinterpolation coefficients for three entire test functions.}
\label{hyper3fig}
\end{figure}

\begin{figure}
\centering
\includegraphics[scale=0.60,clip]{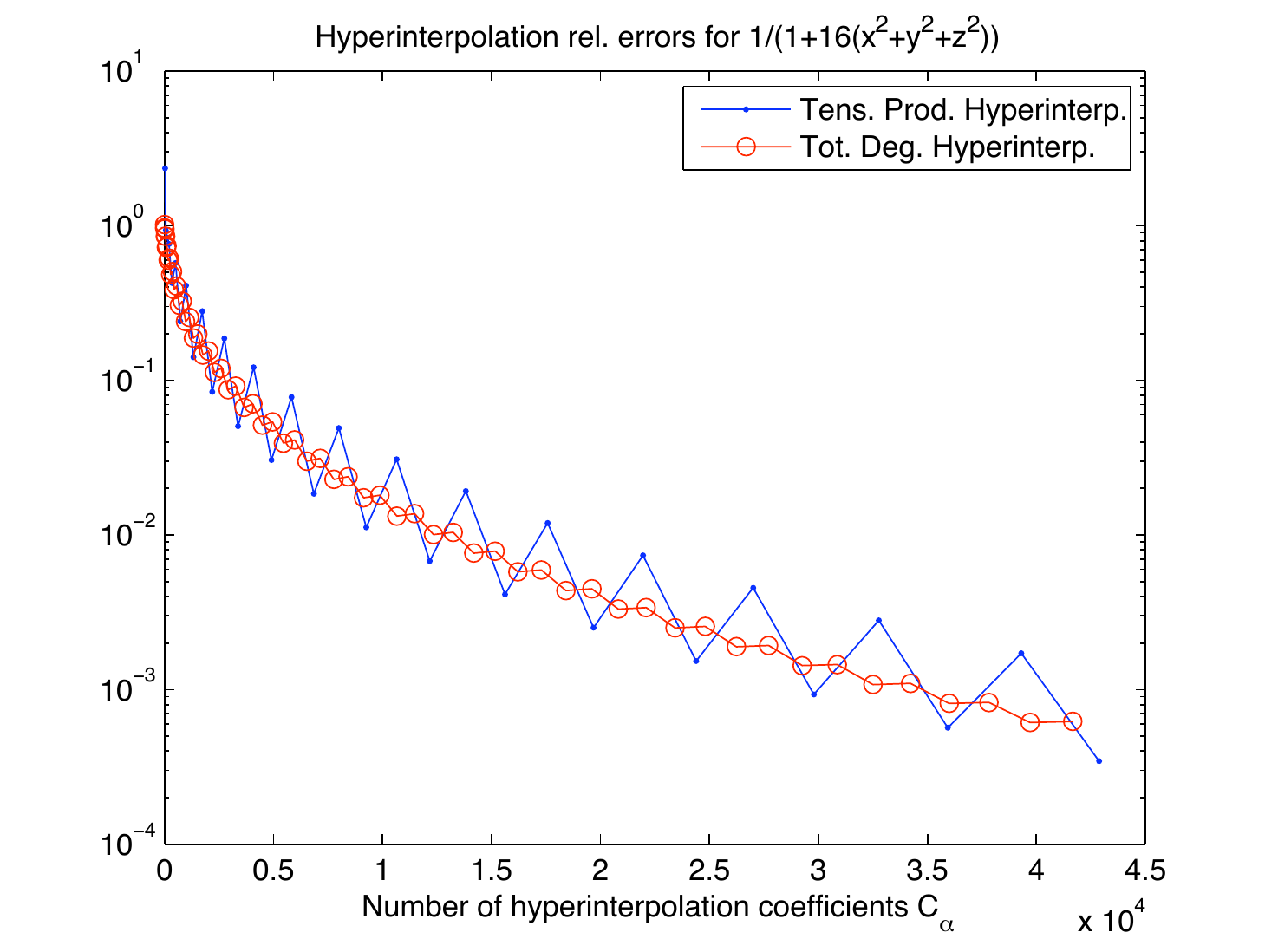}\hfill
\includegraphics[scale=0.60,clip]{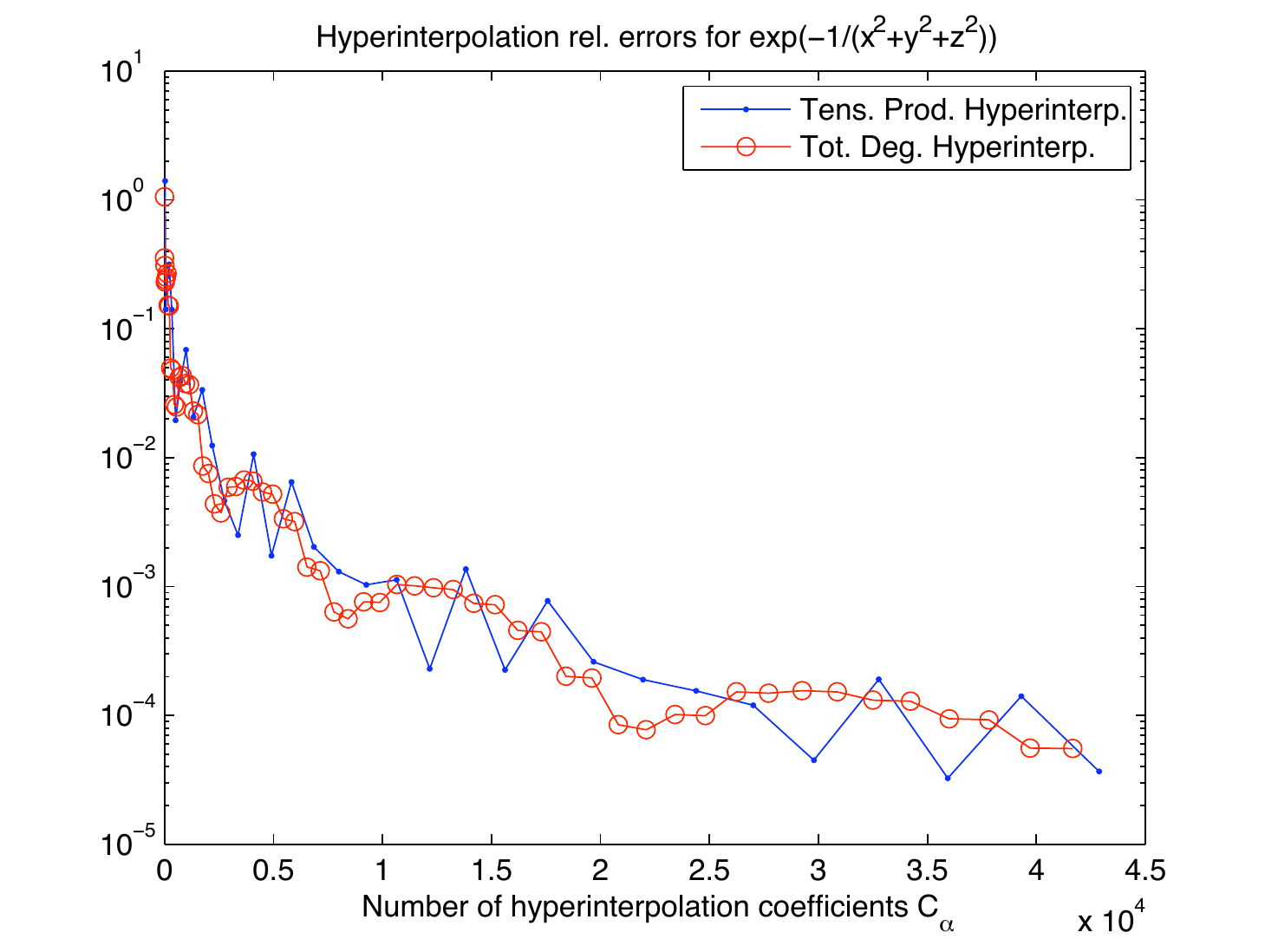}\hfill
\includegraphics[scale=0.60,clip]{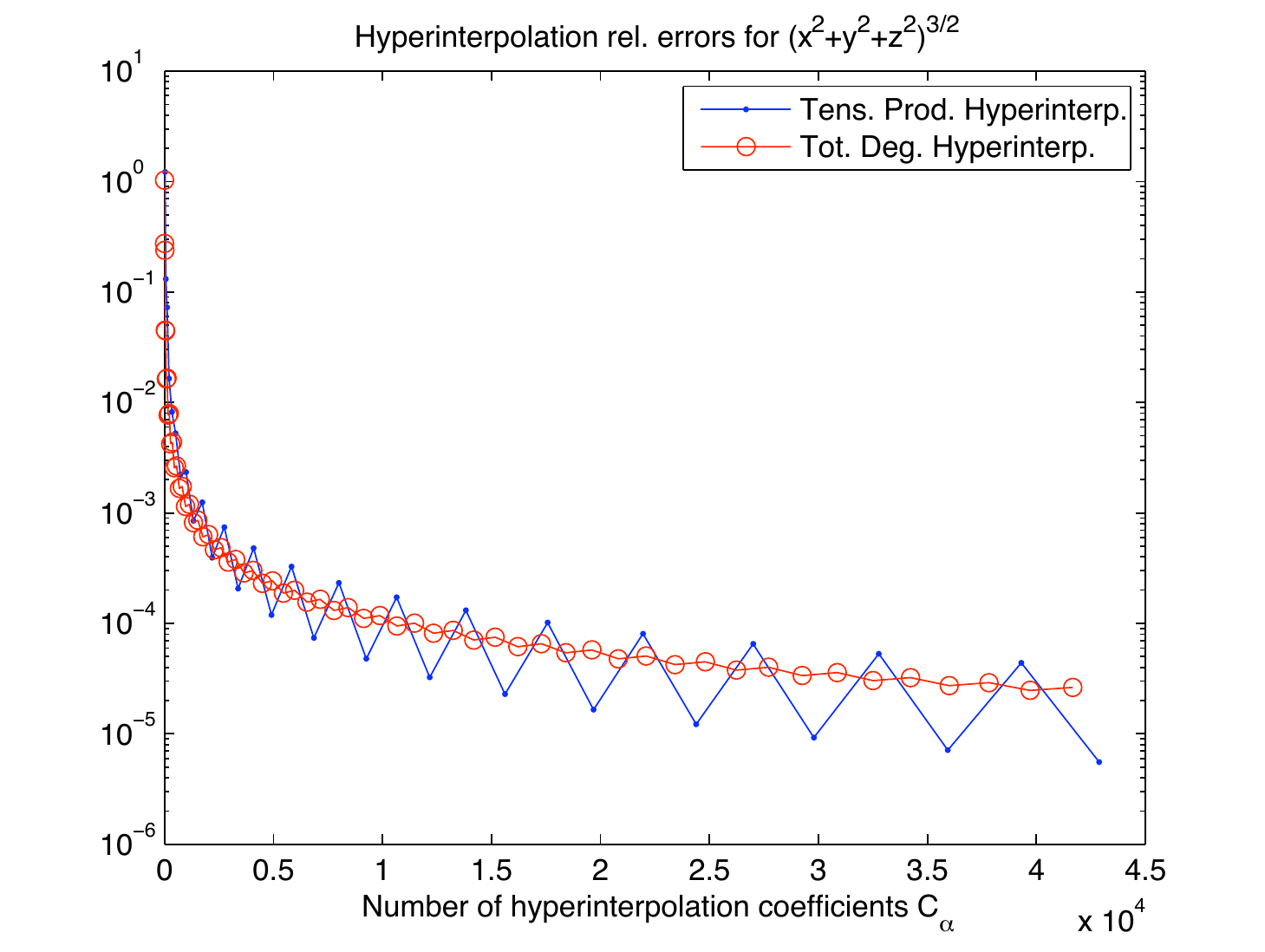}\hfill
\caption{Hyperinterpolation relative errors versus the number of 
hyperinterpolation 
coefficients for 
three nonentire test functions.}
\label{hyper4fig}
\end{figure}

\section{A Clenshaw-Curtis-like formula in the cube.}
In the recent paper \cite{SVZ}, perusing an idea already present 
in \cite{S95}, it has been shown how hyperinterpolation allows us 
to construct new cubature formulae. Given $h\in
L^2_{d\mu}(\Omega)$ and $f\in C(\Omega)$, we can approximate 
the integral of $hf$ in $d\mu$ as
\begin{align} \label{prodcub}
\int_\Omega{h(x)\,f(x)\,d\mu} & \approx
\int_\Omega{h(x)\,\mathcal{L}_n f(x)\,d\mu} \notag \\
&=\sum_{|\alpha|\leq n}{c_\alpha\,m_\alpha}=
\sum_{{\xi}\in X_n}{\lambda_{\xi}\,f({\xi})}\;,
\end{align}
where the generalized ``orthogonal moments'' $\{m_\alpha\}$ and the
cubature weights $\{\lambda_{\xi}\}$ are defined by
\begin{equation} \label{momwei}
m_\alpha:=\int_\Omega{h(x)\,p_\alpha(x)\,d\mu}\;,\;\;
\lambda_{\xi}:=w_{\xi}\sum_{|\alpha|\leq n}{p_\alpha(\xi)\,m_\alpha}\;.
\end{equation}
Observe that the cubature formula (\ref{prodcub}) is {\em exact\/}
for every $f\in \Pi_n^d(\Omega)$, and that $\{m_\alpha\}$ are just
Fourier coefficients of $h$ with respect to the $\mu$-orthonormal basis
$\{p_\alpha\}$. 

Concerning stability and convergence of such cubature formulae, the 
following 
result has been proved in \cite{SVZ}:
\begin{thm}
Let all the assumptions for the construction of the cubature formula
(\ref{prodcub}) be satisfied, and in particular let $h\in
L^2_{d\mu}(\Omega)$.
Then the sum of the absolute values
of the cubature weights has a finite limit
\begin{equation} \label{lim}
\lim_{n\to \infty}{\sum_{{\xi}\in
X_n}{\left |\lambda_{\xi}\right |}}=\int_\Omega{|h(x)|\, d\mu}\;.
\end{equation}
\end{thm} 
\vskip0.3cm

Notice that (\ref{lim}) ensures that the sum of absolute values of
the weights is bounded, and thus by recalling that $\mathcal{L}_n$ 
is a projection operator on $\Pi^d_n(\Omega)$ we obtain the 
Polya-Steklov type (cf. \cite{Kry})  
convergence estimate
\begin{equation} \label{PS}
\left|\int_\Omega{h(x)\,f(x)\,d\mu}-\sum_{{\xi}\in 
X_n}{\lambda_{\xi}\,f({\xi})}\right|\leq \left(\int_\Omega{|h(x)|\, d\mu}
+\sup_n{\sum_{{\xi}\in
X_n}{\left |\lambda_{\xi}\right |}}\right)\,E_n(f)\;,
\end{equation}
where $E_n(f)$ denotes the error of the best polynomial approximation 
of degree $n$ to $f$ in 
the uniform norm.  

Now, applying (\ref{prodcub})-(\ref{momwei}) in the case
\begin{equation} 
d\mu=w(x)\,dx\;,\;\;w\in 
L^1_{dx}(\Omega)\;,\;\;\mbox{with}\;\;h=\frac{1}{w}\in 
L^1_{dx}(\Omega)\;, 
\end{equation}
(since then $h^2=1/w^2\in L^1_{d\mu}(\Omega)$) 
we obtain, via hyperinterpolation, a cubature formula for the standard 
Lebesgue measure from an algebraic cubature 
formula for another measure (absolutely continuos with respect to the 
former). The specialization of this approach to the $1$-dimensional 
Chebyshev measure gives ultimately the popular Clenshaw-Curtis quadrature 
formula \cite{CC}. An extension to dimension 2 has been studied in 
\cite{SVZ}. 
Here we apply the method in dimension 3, obtaining a new nontensorial 
Clenshaw-Curtis-like cubature formula in the 3-cube. 

In Figures 7-8 we display the relative errors of such a formula 
for $(\a_1,\a_2,\a_3)=(E,E,E)$ (cf. (\ref{hypernodes})) on the six 
test functions already used above, compared with 
those of the tensor-product Clenshaw-Curtis, Gauss-Legendre, 
and Gauss-Legendre-Lobatto formulae. The numerical results have been 
obtained with Matlab, using \cite{G04} for the Gaussian formulae 
and \cite{vW} for the tensor-product Clenshaw-Curtis formula.  

In particular, we see that 
with the entire test functions nontensorial Clenshaw-Curtis 
cubature is more accurate than the tensor-product version, 
but less accurate than the other two tensor-product formulae.     
On the other hand, in the less smooth cases the nontensorial 
Clenshaw-Curtis formula is better than all the other three, 
especially for {\em odd\/} hyperinterpolation degrees $n$, which 
correspond to use $n+1$ {\em even\/} in (\ref{cuba-d}) (again a 
sort of parity phenomenon, cf. Figure 2). This behavior echos that of 
1-dimensional and 2-dimensional Clenshaw-Curtis formulae 
(see \cite{T08,SVZ}). Other numerical tests (not reported for brevity) 
have shown that the other versions of the 
nontensorial Clenshaw-Curtis formula, corresponding to 
$(\sigma_1,\sigma_2,\sigma_3)= 
(E,E,O),\,(E,O,E),\,(O,E,E)$ in (\ref{hypernodes}), produce essentially 
the same results.   

\begin{figure}
\centering
\includegraphics[scale=0.60,clip]{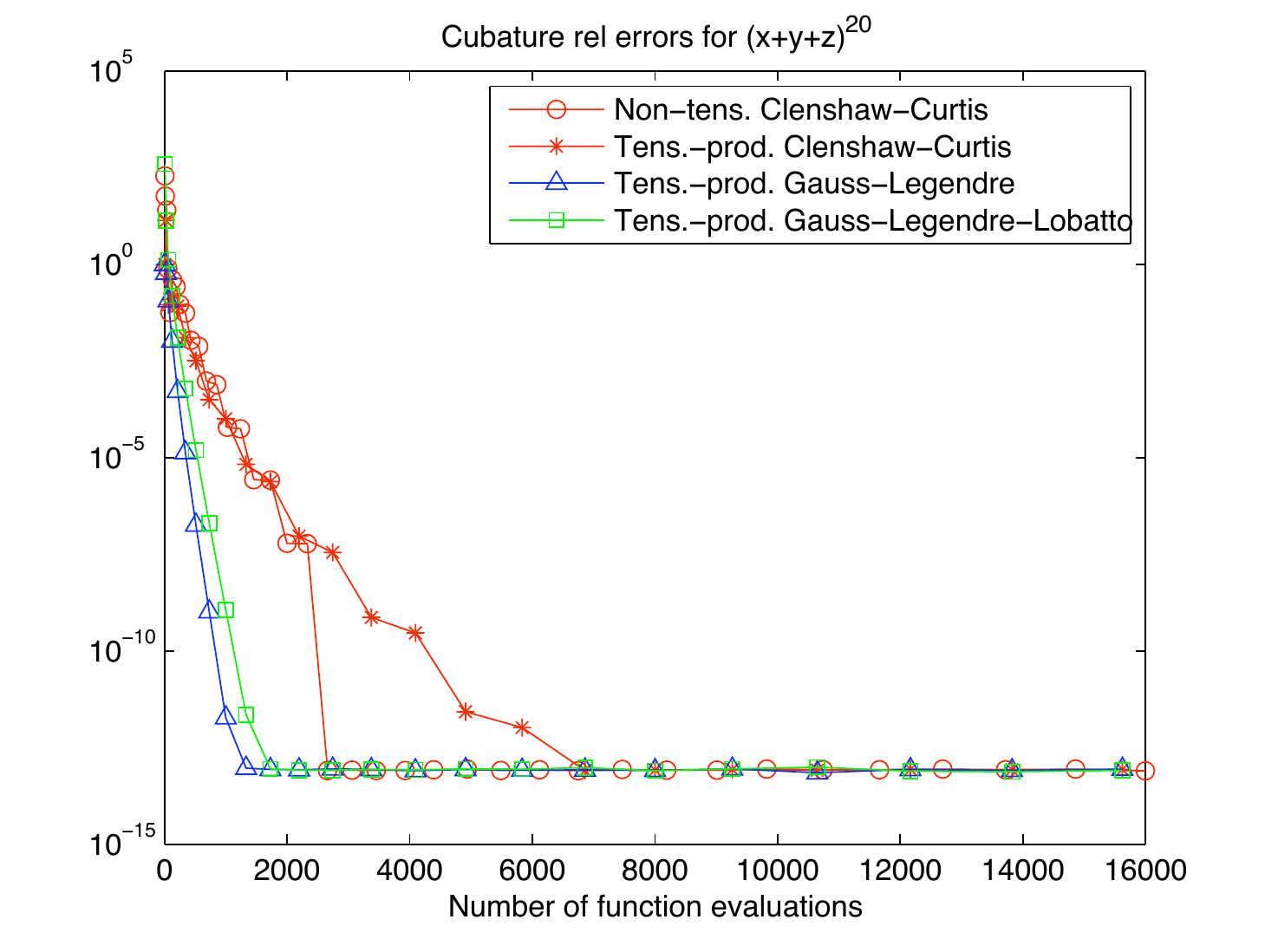}\hfill
\includegraphics[scale=0.60,clip]{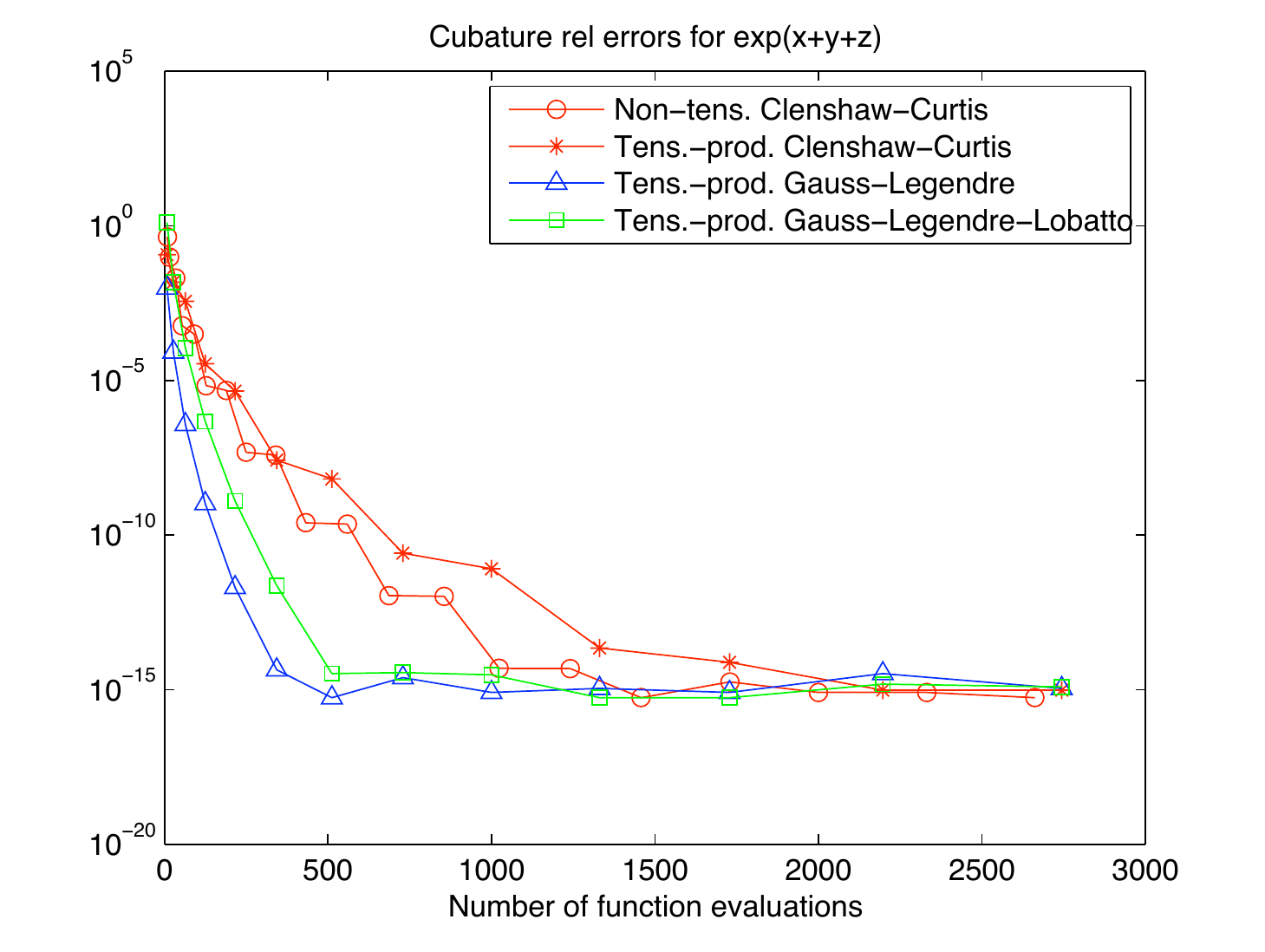}\hfill
\includegraphics[scale=0.60,clip]{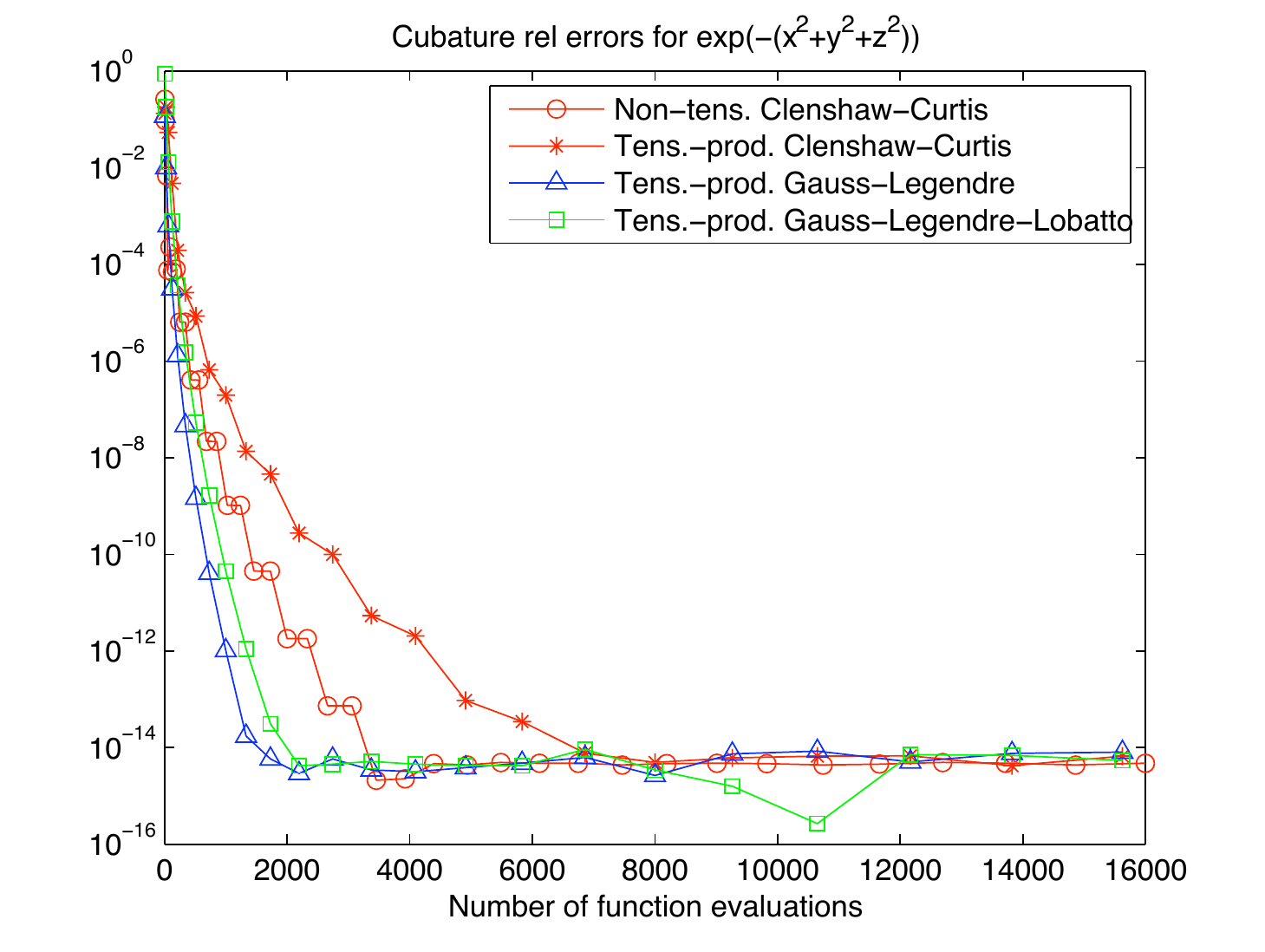}\hfill
\caption{Relative cubature errors versus the number of 
cubature points for three test functions.}
\label{cubCCfig1}
\end{figure}

\begin{figure}
\centering
\includegraphics[scale=0.60,clip]{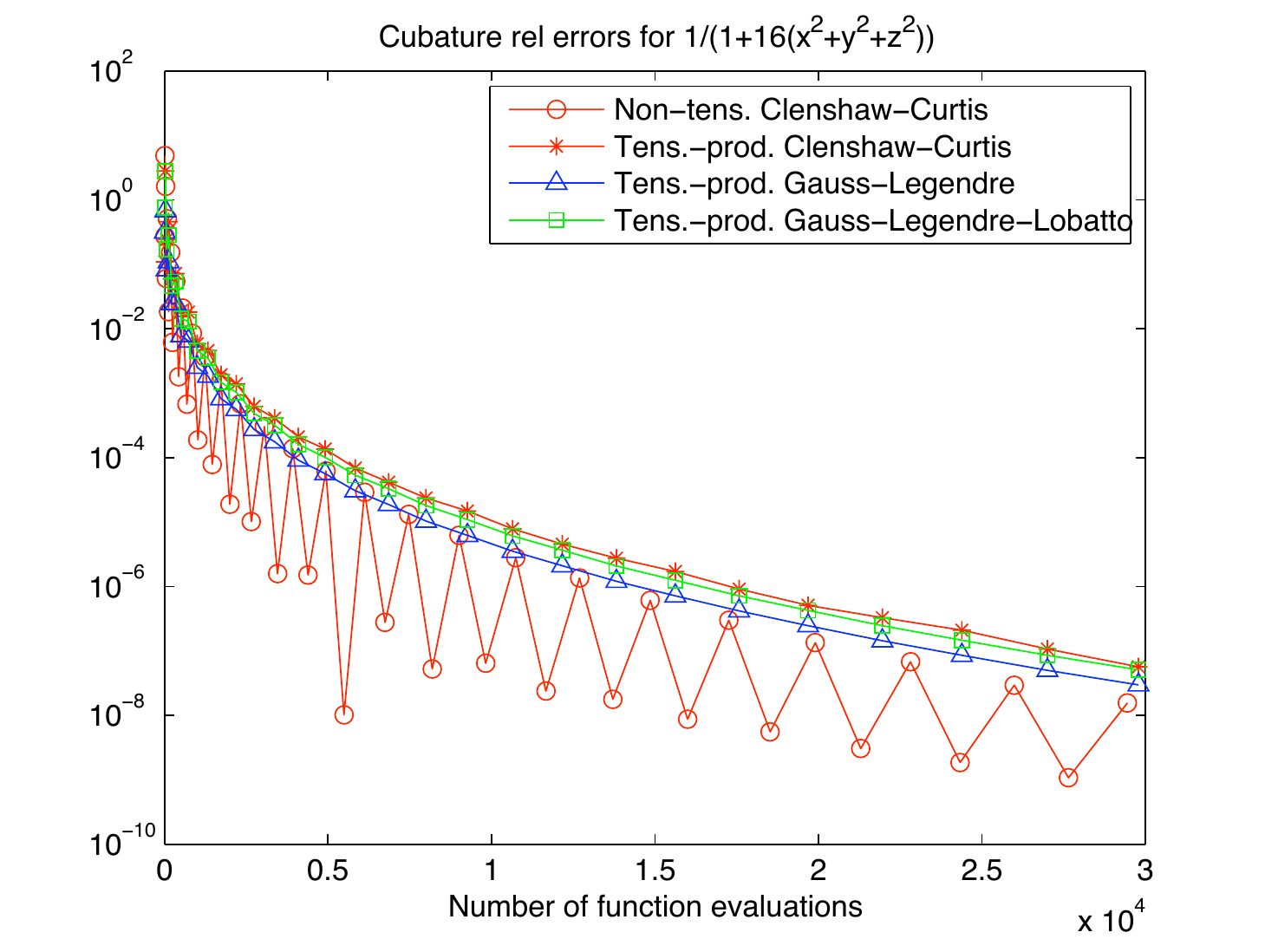}\hfill
\includegraphics[scale=0.60,clip]{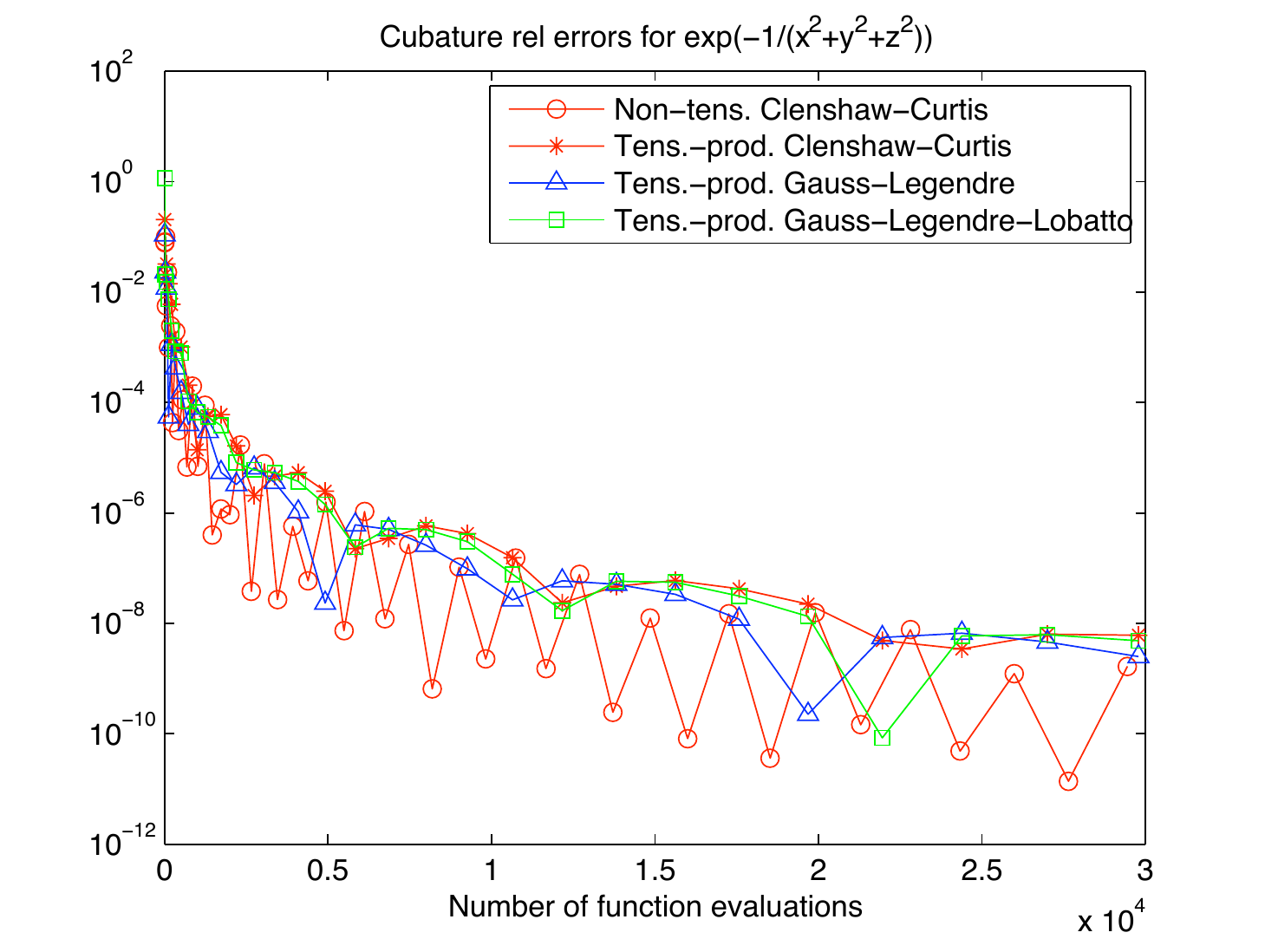}\hfill
\includegraphics[scale=0.60,clip]{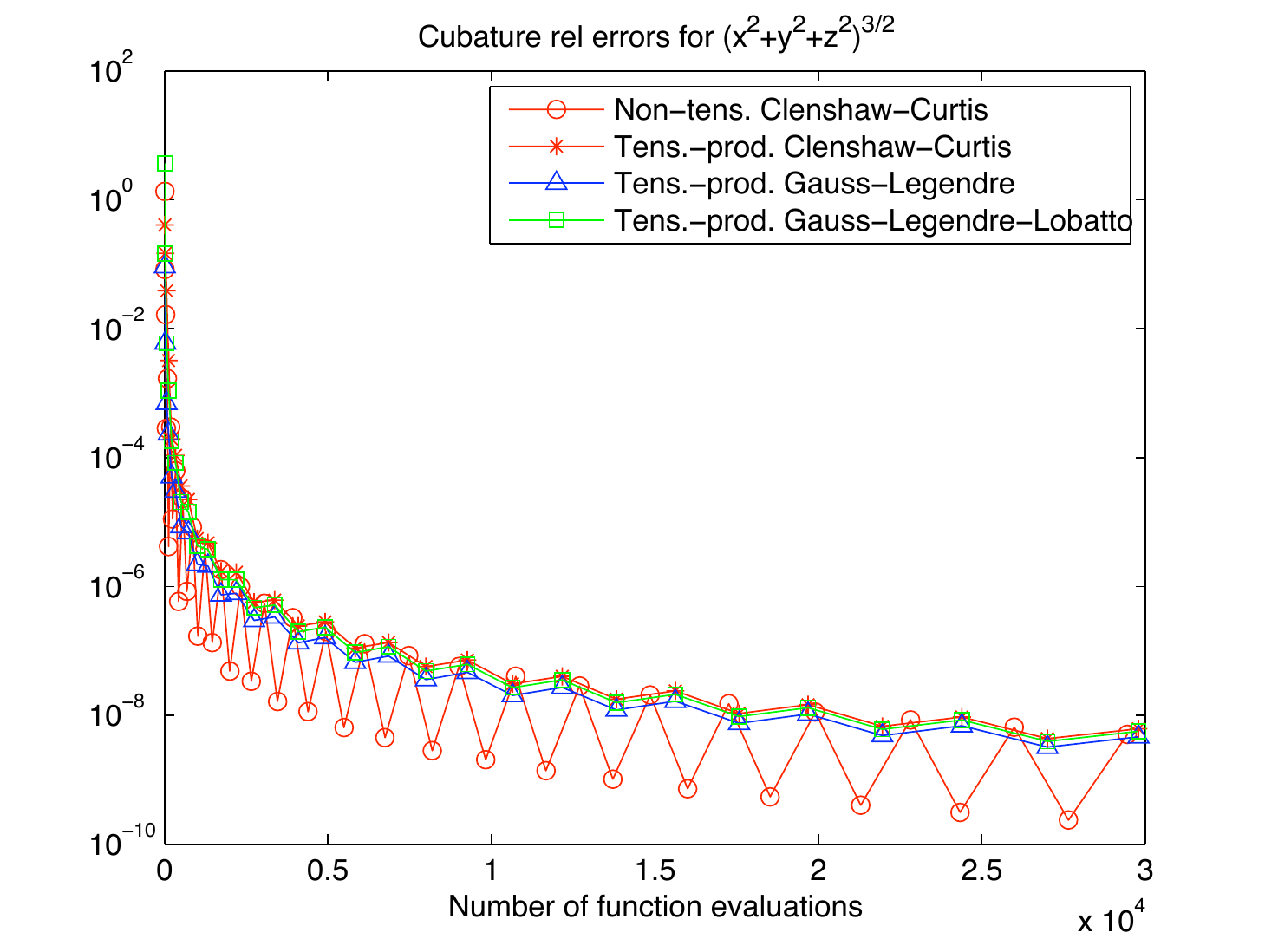}\hfill
\caption{Relative cubature errors versus the number of
cubature points for three test functions.}
\label{cubCCfig2}
\end{figure}

\end{document}